\documentclass[a4paper,11pt]{article}
\usepackage[utf8]{inputenc}
\usepackage{enumerate}
\usepackage{lmodern}
\usepackage[T1]{fontenc}
\usepackage{verbatim}
\usepackage{textcomp}

\usepackage[english]{babel}
\usepackage[a4paper,vmargin={3.5cm,3.5cm},hmargin={2.5cm,2.5cm}]{geometry}
\usepackage[font=sf, labelfont={sf,bf}, margin=1cm]{caption}
\usepackage[pdftex]{hyperref}

\usepackage[pdftex]{color,graphicx}

\usepackage{amsmath,amsfonts,amssymb,amsthm,mathrsfs}

\newcommand{\e}{\mbox{e}\hspace{-.34em}\mbox{\i}}

\newtheorem{theorem}{Theorem}

\newtheorem{lemma}[theorem]{Lemma}

\newtheorem{corollary}[theorem]{Corollary}
\newtheorem{proposition}[theorem]{Proposition}

\def\llbracket{[\hspace{-.10em} [ }
\def\rrbracket{ ] \hspace{-.10em}]}

\def\ww{\mathcal{W}}
\def\w{\mathrm{w}}
\def\ll{{\mathcal L}}
\def\t{{\mathcal T}}
\def\cc{\mathcal{C}}
\def\ve{\varepsilon}
\def\wt{\widetilde}
\def\wh{\widehat}
\def\bm{{\bf m}}

\def\bn{{\bf n}}
\def\bp{{\bf p}}
\def\ff{\mathcal{F}}
\def\nn{\mathcal{N}}
\def\mm{\mathcal{M}}
\def\dd{\mathcal{D}}
\def\yy{\mathcal{Y}}
\def\zz{\mathcal{Z}}
\def\bN{{\bf N}}
\def\D{\mathrm{d}}
\def\z{\mathcal{Z}}
\def\e{\mathcal{E}}
\def\la{\longrightarrow}
\def\E{{\mathbb E}}
\def\P{{\mathbb P}}
\def\N{{\mathbb N}}
\def\R{{\mathbb R}}
\def\Q{{\mathbb Q}}

\def\ov{\overline}
\def\un{\underline}
\def\build#1_#2^#3{\mathrel{
\mathop{\kern 0pt#1}\limits_{#2}^{#3}}}
\def\rem{\noindent{\bf Remark. }}

\title{Subordination of trees and the Brownian map}
\author{Jean-Fran\c cois Le Gall}

\date{\small \it Universit\'e Paris-Sud}
\begin{document}
\maketitle

\begin{abstract}
We discuss subordination of random compact $\R$-trees. We focus on
the case of the Brownian tree, where the subordination function 
is given by the past maximum process of Brownian motion indexed
by the tree. In that particular case, the subordinate tree is 
identified as a stable L\'evy tree with index $3/2$. As a more precise
alternative formulation, we show that the maximum process of the Brownian snake
is a time change of the height process coding the L\'evy tree.
We then apply our results to properties of the Brownian map.
In particular, we recover, in a more precise form, a recent result 
of Miller and Sheffield identifying the metric net associated with
the Brownian map.
\end{abstract}

\section{Introduction}

Subordination is a powerful tool in the study of random processes.
In the present work, we investigate subordination of random trees, and we 
apply our results to properties of the random metric space called the Brownian map,
which has been proved to be the universal scaling limit of many different
classes of random planar maps (see in particular \cite{AA,Uniqueness,Mie}). These applications have been motivated
by the work of Miller and Sheffield \cite{MS}, which is part of a program aiming at the construction of
a conformal structure on the Brownian map (see \cite{MS1,MS2} for recent
developments in this direction).

\begin{figure}[!h]
 \begin{center}
 \includegraphics[width=14cm]{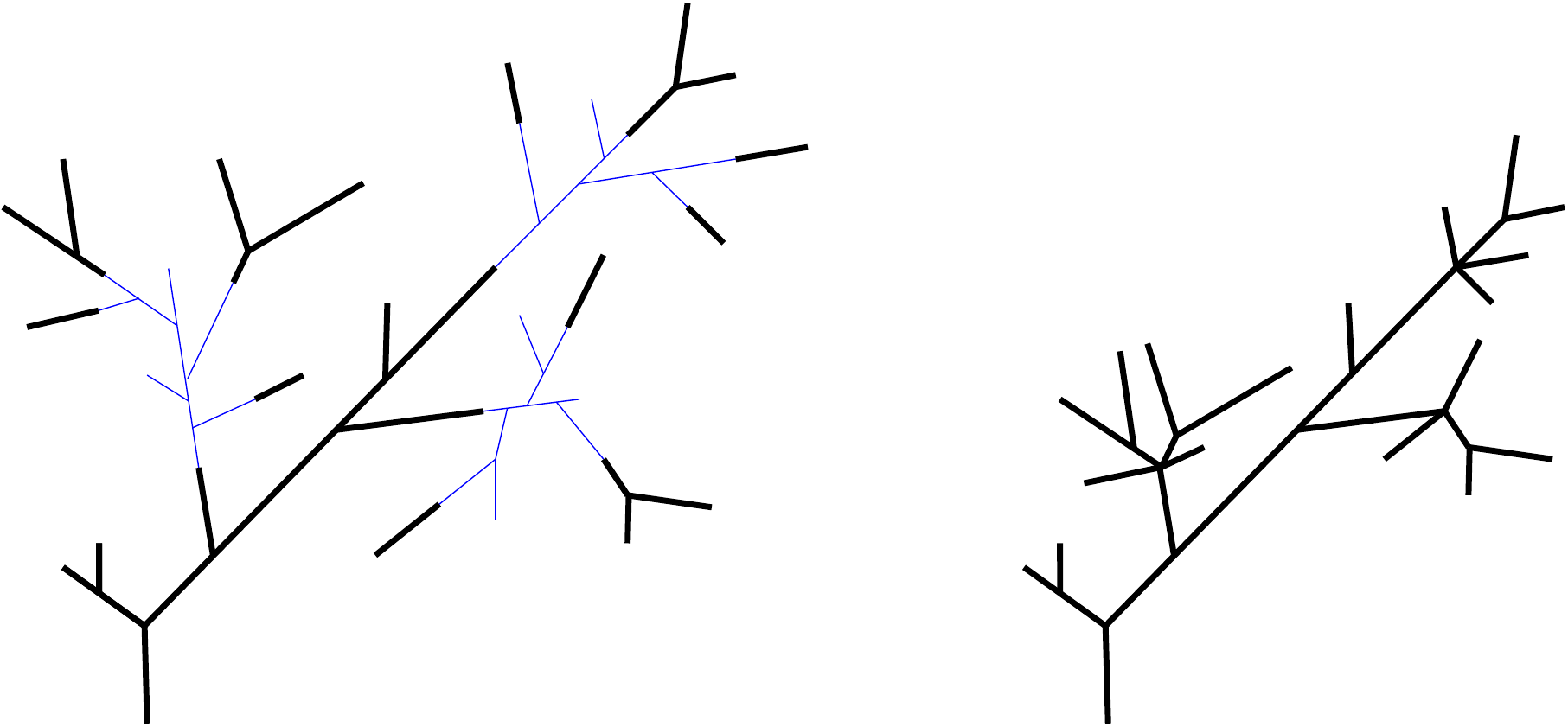}
 \caption{\label{subordina-tree}
On the left side, the tree $\t$, with segments where $g$ is constant pictured in thin blue lines.
On the right side, the subordinate tree, where each connected component of $\t$ made of
thin segments has been glued into a single point.}
 \end{center}
 \vspace{-8mm}
 \end{figure}

To explain our starting point, let us consider a compact $\R$-tree $\t$. This means 
that $\t$ is a compact metric space such that, for every $a,b\in\t$,
there exists a unique (continuous injective) path from $a$ to $b$,
up to reparameterization, and the range of this path, which is called
the geodesic segment between $a$ and $b$ and denoted by $\llbracket a,b \rrbracket$,
is isometric to a compact interval of the real line. We assume that $\t$
is rooted, so that there is a distinguished point $\rho$ in $\t$. This allows
us to define a generalogical order on $\t$, by saying
that $a\prec b$ if and only if $a\in \llbracket \rho,b\rrbracket$. Consider
then a continuous function $g:\t\la\R_+$, such that $g(\rho)=0$
and $g$ is nondecreasing for the genealogical order. The basic idea
of subordination is to identify $a$ and $b$ if $g$ is constant on 
the geodesic segment $\llbracket a,b \rrbracket$. So, for every
$a\in\t$, the set of all points that are identified with $a$ is a closed connected subset of
$\t$. This glueing operation yields another compact $\R$-tree $\wt\t$, which is equipped
with a metric such that the distance between $\rho$ and $a$ is $g(a)$ and
is called the subordinate tree of $\t$ with respect to $g$ (see Fig.\ref{subordina-tree} for an 
illustration). Furthermore,
if our initial tree $\t$ was given as the tree coded by a continuous function 
$h:[0,\sigma]\la\R_+$ (see \cite{DLG} or Section \ref{sec:determin} below), the 
subordinate tree $\wt\t$ is coded by $g\circ p_h$, where
$p_h$ is the canonical projection from $[0,\sigma]$ onto $\t$.

Our main interest is in the case where $\t$ is random, and more 
precisely $\t=\t_\zeta$ is the ``Brownian tree'' coded by a positive Brownian 
excursion $\zeta=(\zeta_s)_{0\leq s\leq \sigma}$ under the It\^o excursion measure.
One may view $\t_\zeta$ as a variant of Aldous' Brownian CRT, for which the total mass is
not finite, but is distributed according to an infinite measure on $(0,\infty)$. 
As previously, we write  $p_\zeta$ for the canonical projection from $[0,\sigma]$ onto $\t_\zeta$.
Next, to define the subordination function, we let $(Z_a)_{a\in \t_\zeta}$
be (linear) Brownian motion indexed by $\t_\zeta$, starting from $0$ at the root $\rho$.
A simple way to construct this process is to use the Brownian snake approach, letting 
$Z_a=\wh W_s$ if $a=p_\zeta(s)$, where $\wh W_s$ is the ``tip'' of the random path $W_s$,
which is the value at time $s$ of the Brownian snake driven by $\zeta$. Since $\zeta$
is distributed according to the It\^o measure, $W$ follows the Brownian snake excursion 
measure away from $0$, which we denote by $\N_0$ (see \cite[Chapter IV]{Zurich}). We then set 
$\ov Z_a=\max\{Z_b : b\in\llbracket \rho,a \rrbracket\}$. In terms
of the Brownian snake, we have $\ov Z_a=\ov W_s:=\max\{ W_s(t):0\leq t\leq \zeta_s\}$ 
whenever $a=p_\zeta(s)$. We also use the notation $\un W_s:=\min\{ W_s(t):0\leq t\leq \zeta_s\}$.

\begin{theorem}
\label{subor-brown}
Let $\wt\t_\zeta$ stand for the 
subordinate tree of $\t_\zeta$ with respect to the (continuous nondecreasing) function $a\mapsto \ov Z_a$.
Under the Brownian snake excursion measure $\N_0$, the tree $\wt \t_\zeta$ is a L\'evy tree with branching mechanism
$$\psi_0(r)=\sqrt{\frac{8}{3}} \, r^{3/2}.$$
\end{theorem}

Recall that L\'evy trees represent the genealogy of continuous-state branching
processes \cite{DLG}, and can be characterized by a regenerative property analogous to
the branching property of Galton--Watson trees \cite{Weill}. Our identification of the
distribution of $\wt\t_\zeta$ is reminiscent of the classical result stating that
the right-continuous inverse of the maximum process of a standard linear Brownian motion
is a stable subordinator with index $1/2$. 

In view of our applications, it turns out that it is important to have more information than the 
mere identification of $\wt\t_\zeta$ as a random compact $\R$-tree. As mentioned above,
$\wt\t_\zeta$ can be viewed as the tree coded by the random function $s\mapsto 
\ov Z_{p_\zeta(s)}=\ov W_s$.
This coding induces a ``lexicographical'' order structure on $\t_\zeta$ (see \cite{Duq} for
a thorough discussion of order structures on $\R$-trees). Somewhat surprisingly, it 
is not immediately clear that the order structure on $\wt\t_\zeta$ 
induced by the coding from the function $s\mapsto\ov W_s$ coincides with the usual
order structure on L\'evy trees, corresponding to the uniform random
shuffling at every node in the terminology of \cite{Duq}. To obtain this property,
we relate the function $s\mapsto\ov W_s$ to the height process of  \cite{DLG0,DLG}.
We recall that, from a L\'evy process with Laplace exponent $\psi_0$, we can construct
a continuous random process $(H_t)_{t\geq0}$ called the height process, which codes the $\psi_0$-L\'evy tree
(here and below we say $\psi_0$-L\'evy tree rather than L\'evy tree with branching mechanism $\psi_0$ for simplicity). See Section \ref{sec:Levy-tree} below for more details.

\begin{theorem}
\label{subordi-coding}
There exists a process $H$, which is distributed under $\N_0$ as the height process of
a L\'evy tree with branching mechanism $\psi_0(r)=\sqrt{8/3}\,r^{3/2}$, and a continuous random process
$\Gamma$ with nondecreasing sample paths, such that we have, $\N_0$ a.e. for
every $s\geq 0$,
$$
\ov W_s= H_{\Gamma_s}.$$
\end{theorem}

Both $H$ and $\Gamma$ will be constructed  in the proof
of Theorem \ref{subordi-coding}, and
are measurable functions of $(W_s)_{s\geq 0}$. 
It is possible to identify the random process $\Gamma$ as a continuous additive functional 
of the Brownian snake, but we do not need this fact in the subsequent applications, and
we do not discuss this matter in the present work.

\smallskip

Theorem \ref{subordi-coding} implies that the tree coded by $\ov W_s$ is isometric to the tree coded by $H_s$, and we recover 
Theorem \ref{subor-brown}. But 
Theorem \ref{subordi-coding} gives much more, namely that the order structure induced 
by the coding via $s\to \ov W_s$ is the same as the order structure induced by the 
usual height function of the L\'evy tree.
Order structures are crucial for our applications to the Brownian map.
Similarly as in \cite{MS}, we deal with a version of the Brownian map with randomized volume,
which is constructed as a quotient space of the tree $\t_\zeta$: Two points 
$a$ and $b$ of $\t_\zeta$ are identified if $Z_a=Z_b$ and if $Z_c\geq Z_a$
for every $c\in [a,b]$, where $[a,b]$ is the set of all points that are visited when 
going from $a$ to $b$ around the tree in ``clockwise'' order (for this to make sense,
it is essential that $\t_\zeta$ has been equipped with a lexicographical order structure).
We write $\bm$ for the resulting quotient space (the Brownian map) and 
$D^*$
for the metric on the Brownian map (see \cite{Geo,Uniqueness} for more details).
The space $\bm$ comes with two distinguished points, namely the root $\rho$
of $\t_\zeta$ and the unique point $\rho_*$ where $Z$ attains its minimal
value -- in a sense that can be made precise, these two points are
independent and uniformly distributed over $\bm$. Let $B(r)$ be the closed ball of radius $r$ centered at $\rho_*$ in $\bm$. For every
$r\in (0,D(\rho^*,\rho))$, define the {\it hull} $B^\bullet(r)$ as the complement
of the connected component of the complement of $B(r)$ that contains $\rho$ (informally, $B^\bullet(r)$ is
obtained by filling in all holes of $B(r)$ except for the one containing $\rho$). Following 
\cite{MS}, we define the metric net $\mm$ as the closure of
$$\bigcup_{0\leq r<D(\rho^*,\rho)} \partial B^\bullet(r).$$
(This definition is in fact a little different from \cite{MS} which does not take the closure of the union
in the last display.) We can equip the set $\mathcal{M}$ with an ``intrinsic'' metric $\Delta^*$
 derived from 
the Brownian map metric $D^*$.

\begin{figure}[!h]
 \begin{center}
 \includegraphics[width=12cm]{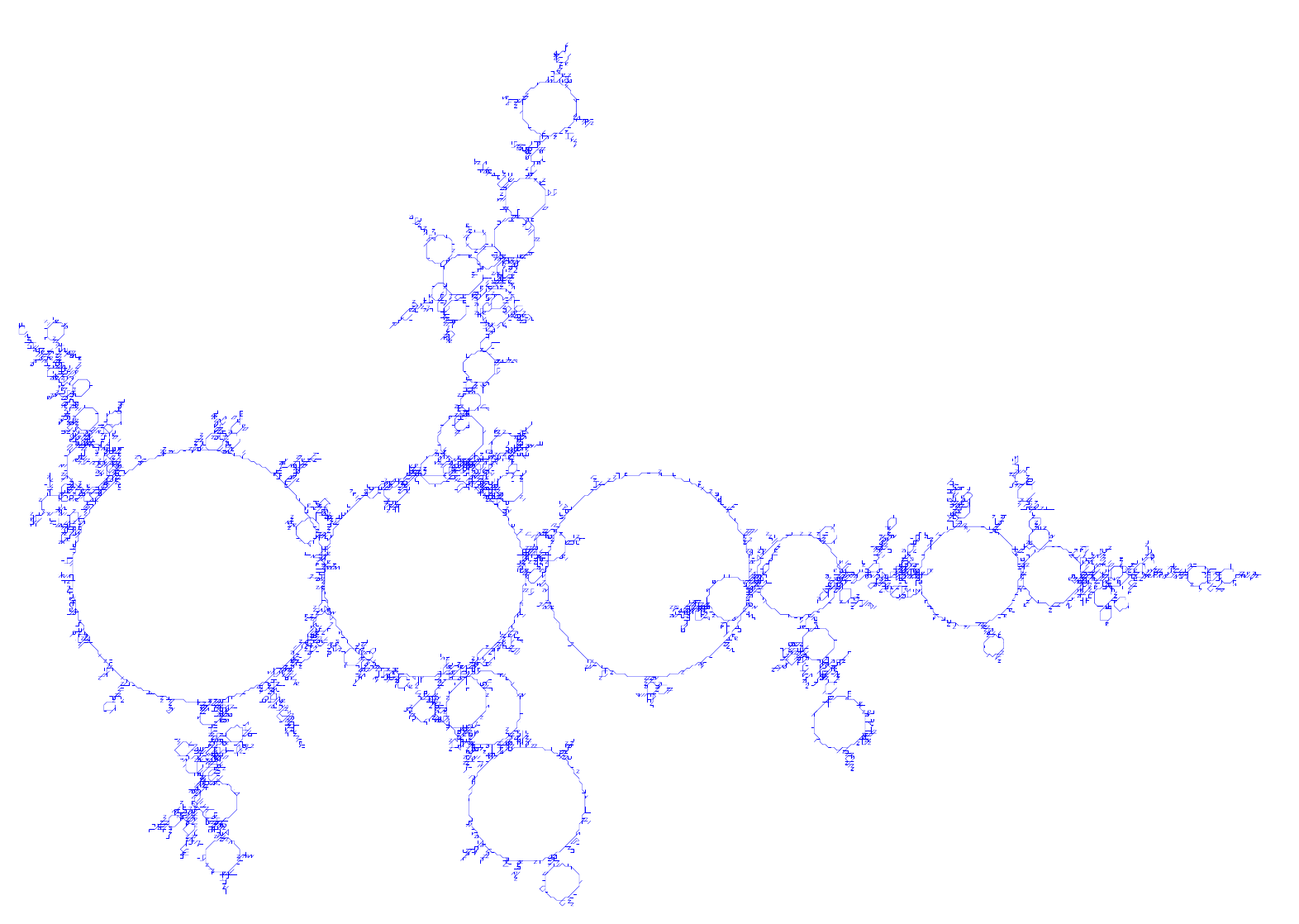}
 \caption{\label{looptree}
 A simulation of a looptree (simulation by Igor Kortchemski). For technical reasons, some of the
 trees branching off a loop are pictured inside this loop, but, from the point of view
 of the present work, it is better to think of these trees as growing outside the loop, so 
 that the space inside the loop may be ``filled in'' appropriately.}
 \end{center}
 \vspace{-6mm}
 \end{figure}

It is not hard to verify that, in the construction of $\bm$ as
a quotient space of $\t_\zeta$, points of $\mm$
exactly correspond to vertices $a$ of $\t_\zeta$
such that $Z_a=\un Z_a:=\min\{Z_b:b\in \llbracket \rho,a\rrbracket\}$
(Proposition \ref{descri-net}).
This suggests that the metric net $\mm$ is closely related to the 
subordinate tree of $\t_\zeta$ with respect to the function
$a\mapsto -\un Z_a$, which is a $\psi_0$-L\'evy tree by the preceding results. 
To make this relation precise, we need the notion of the looptree introduced by
Curien and Kortchemski \cite{CK} in a more general setting.  Informally, the 
looptree associated with a L\'evy tree is constructed by
replacing each point $a$ of infinite multiplicity 
by a loop of ``length'' equal to the ``weight'' of $a$, in such a way that the subtrees that are
the
connected components of the complement of $a$  branch along this loop in
an order determined by the coding function (see Fig.\ref{looptree} below for an
illustration). To give a more precise definition, 
first note that Theorem \ref{subordi-coding}
and an obvious symmetry argument allow us to find a process
$(H'_s)_{s\geq 0}$ distributed as the height process of a
$\psi_0$-L\'evy tree, and 
a continuous random process $(\Gamma'_s)_{s\geq 0}$
with nondecreasing sample paths such that $-\un W_s= H'_{\Gamma'_s}$.  Then let $\chi:=\sup\{s\geq 0:H'_s>0\}$, and let $\sim$
be the equivalence relation on $[0,\chi]$ whose graph is the smallest closed symmetric subset 
of $[0,\chi]^2$ that contains all pairs $(s,t)$ with $s\leq t$, $H'_s=H'_t$ and
$H'_r>H'_s$ for every $r\in (s,t)$. The looptree $\ll$ is defined as the quotient space
$[0,\chi]/\sim$ equipped with an appropriate metric (see \cite{CK}). Clearly,
$H'_\alpha$ makes sense for any $\alpha\in \ll$ and is interpreted as the height of
$\alpha$. The metric on $\ll$
is in fact not relevant for us, since we consider instead the pseudo-metric $\dd^\circ$
defined for $\alpha,\beta\in\ll$ by
$$\dd^\circ(\alpha,\beta)= 2\min\Bigg(
\max_{\gamma\in[\alpha,\beta]} H'_\gamma, 
\max_{\gamma\in[\beta,\alpha]} H'_\gamma\Bigg) - H'_\alpha -H'_\beta,$$
where $[\alpha,\beta]$ corresponds to the subset of $\ll$ visited when going from
$\alpha$ to $\beta$ ``around'' $\ll$ in ``clockwise order'' (see Section \ref{sec:Br-map}
for a more formal definition). We write $\alpha\simeq \beta$ if $\dd^\circ(\alpha,\beta)=0$
(informally this means that $\alpha$ and $\beta$ ``face each other'' in the tree,
in the sense that they are at the same 
height, and that points ``between'' $\alpha$ and $\beta$ are at a smaller height).
It turns out that this defines an equivalence relation on $\ll$. Finally, we let
$\dd^*$ be the largest symmetric function on $\ll\times \ll$ that is bounded above by
$\dd^\circ$ and satisfies the triangle inequality.

\begin{theorem}
\label{identi-net}
The metric net $(\mm, \Delta^*)$ is a.s. isometric to the quotient space 
$\ll/\!\simeq$ equipped with the metric induced by $\dd^*$.
\end{theorem}

See Theorem \ref{identi-metricnet} in Section \ref{sec:Br-map}
for a more precise formulation. Theorem \ref{identi-net} is closely related to
Proposition 4.4 in \cite{MS}, where, however, the metric net is not identified as
a metric space. The description of the metric net is an important ingredient
of the axiomatic characterization of the Brownian map discussed in \cite{MS}. 

Let us briefly comment on the motivations for studying the metric net. Roughly speaking, the 
Brownian map $\bm$ can be recovered from the metric net $\mm$
by filling in the ``holes''. To make this precise, we observe that the connected components
of $\bm\backslash \mm$ are bounded by Jordan curves (Proposition \ref{connec-compo})
and that these components are in one-to-one correspondence with connected 
components of $\t_\zeta\backslash \Theta$, where $\Theta:=\{a\in\t_\zeta: Z_a=\un Z_a\}$. 
Each of the latter components is associated with an excursion of the 
Brownian snake above its minimum, in the terminology of \cite{ALG}, and the distribution
of such an excursion only depends on its boundary size as defined in \cite{ALG} (this boundary size
can be interpreted as a generalized length of the Jordan curve bounding the 
corresponding component of $\bm\backslash \mm$). Theorem 40 in \cite{ALG}
shows that, conditionally on their boundary sizes, these excursions are independent 
and distributed according to a certain ``excursion measure''. In the Brownian map
setting, this means that the holes in the metric net are filled in independently,
conditionally on the lengths of their boundaries.

The paper is organized as follows. Section 2 gives a brief discussion of subordination for deterministic trees, and
Section 3 recalls the basic facts about L\'evy trees that we need.
After a short presentation of the Brownian snake, Section 4 gives the distribution
of the subordinate tree $\wt\t_\zeta$ (Theorem \ref{subor-brown}). In view of identifying the order structure
of this subordinate tree, Section 5 provides a technical result showing that the 
height process coding a L\'evy tree is the limit in a strong sense of the
discrete height functions coding embedded Galton--Watson trees. This result
is related to the general limit theorems of \cite[Chapter 2]{DLG0} proving that 
L\'evy trees are weak limits of Galton--Watson trees, but the fact that we get a strong 
approximation is crucial for our applications. In Section 6, we prove that the 
Brownian snake maximum process $s\mapsto \ov W_s$
is a time change of the height process associated with a $\psi_0$-L\'evy tree
(Theorem \ref{subordi-coding}).
This result is a key ingredient of the developments of Section 7, where we identify the metric net of the Brownian map (Theorem \ref{identi-net}). Section 8 discusses the
connected components of the complement of the metric net, showing 
in particular that they are in one-to-one correspondence with the
points of infinite multiplicity of the associated L\'evy tree, and that the boundary 
of each component is a Jordan curve. Section 9, which is
mostly independent of the preceding sections, discusses more general subordinations 
of the Brownian tree $\t_\zeta$, which lead to stable L\'evy trees with an arbirary index. This 
section is related to our previous article \cite{BLGLJ}, which dealt with 
subordination for spatial branching processes, but the latter work
did not consider the associated genealogical structures as we do here, and
the subordination method, based on the so-called residual lifetime process,
was also different. Finally, the appendix presents a more general and more precise version 
of the special Markov property of the Brownian snake (first established in \cite{snakesolutions}),
which plays an important role in several proofs.

\section{Subordination of deterministic trees}
\label{sec:determin}

In this short section, we present a few elementary considerations about deterministic
$\R$-trees. We refer to \cite{Evans} for the basic facts about $\R$-trees that we will need, and
to \cite{Duq} for a thorough study of the coding of compact $\R$-trees by functions.

Let us consider a compact $\R$-tree $(\t,d)$ rooted at $\rho$. If $a,b\in\t$, the 
geodesic segment between $a$ and $b$ (the range of the unique geodesic from 
$a$ to $b$) is denoted by $\llbracket a,b \rrbracket$. The point $a\wedge b$ is then defined by
$\llbracket \rho,a\wedge b \rrbracket=\llbracket \rho, a \rrbracket\cap \llbracket \rho,b \rrbracket$.
The genealogical partial order on $\t$ is denoted by $\prec\;$: we have $a\prec b$ if 
and only if $a\in \llbracket \rho,b\rrbracket$, and we then say that $a$ is an
ancestor of $b$, or $b$ is a descendant of $a$. Finally, the height of $\t$ is defined by
$$\mathcal{H}(\t)=\max_{a\in\t} d(\rho,a).$$

Let $g:\t\la \R_+$ be a nonnegative continuous function on $\t$. Assume that $g(\rho)=0$
and that $g$ is nondecreasing with respect to the genealogical order ($a\prec b$ implies that
$g(a)\leq g(b)$). We then define, for every $a,b\in\t$,
$$d^{(g)}(a,b)= g(a)+g(b)-2 \,g(a\wedge b).$$
Notice that $d^{(g)}(a,b)$ is a symmetric function of $a$ and $b$ and satisfies the
triangle inequality. We can thus consider the equivalence relation
$$a\approx_{g} b\quad\hbox{if and only if}\quad d^{(g)}(a,b)=0.$$
Thus $a\approx_{g} b$ if and only if $g(a)=g(b)=g(a\wedge b)$, and this is also equivalent to
saying that $g(c)=g(a)$ for every $c\in\llbracket a,b\rrbracket$.
Write $\t^{(g)}$ for the quotient $\t/\!\approx_{g}$, and $\pi^{(g)}$ for the canonical projection from 
$\t$ onto $\t^{(g)}$.

\begin{proposition}
\label{tree-quotient}
The set $\t^{(g)}$ equipped with the distance induced by $d^{(g)}$ is again
a compact $\R$-tree.
\end{proposition}

\proof One immediately verifies that, for every $a\in \t$, $\pi^{(g)}(\llbracket \rho,a \rrbracket)$
is a segment in $(\t^{(g)},d^{(g)})$ with endpoints $\pi^{(g)}(\rho)$ and $\pi^{(g)}(a)$. From 
Lemma 3.36 in \cite{Evans}, it then suffices to check that the four-point condition
$$d^{(g)}(a_1,a_2) + d^{(g)}(a_3,a_4)
\leq \max\Big( d^{(g)}(a_1,a_3)+ d^{(g)}(a_2,a_4), d^{(g)}(a_1,a_4)+d^{(g)}(a_2,a_3)\Big)$$
holds for every $a_1,a_2,a_3,a_4\in\t$. This is straightforward and we omit the details. \endproof

We call $(\t^{(g)},d^{(g)})$ the subordinate tree of $(\t,d)$ with respect to the 
function $g$. By convention, $\t^{(g)}$ is rooted at $\pi^{(g)}(\rho)$. Since 
$d^{(g)}(\rho,a)= g(a)$ for every $a\in\t$, we have
$\mathcal{H}(\t^{(g)})=\max\{g(a):a\in\t\}$.

Consider now a continuous function $h:[0,\sigma]\la \R_+$, where $\sigma \geq 0$, 
such that $h(0)=h(\sigma)=0$, and
assume that $(\t,d)$ is the tree coded by $h$ in the sense of \cite{DLG} or \cite{Duq}. This means that
$\t=\t_h$ is the quotient space $[0,\sigma]/\!\sim_h$, where the equivalence 
relation $\sim_h$ is defined on $[0,\sigma]$ by
$$s\sim_h t \quad\hbox{if and only if}\quad h(s)=h(t)=\min_{s\wedge t\leq r\leq r\vee t} h(r),$$
and $d=d_h$ is the distance induced on the quotient space  by
$$d_h(s,t)= h(s)+h(t)-2\min_{s\wedge t\leq r\leq r\vee t} h(r).$$
Notice that the topology of $(\t_h,d_h)$ coincides with the quotient topology on $\t_h$.

\smallskip

The canonical projection from $[0,\sigma]$ onto $\t_h$ is denoted by $p_h$, and
$\t_h$ is rooted at $\rho_h=p_h(0)$. For $s\in[0,\sigma]$, the quantity 
$h(s)=d_h(0,s)$ is interpreted as the height of $p_h(s)$ in the tree. One easily verifies that, for every $s,t\in[0,\sigma]$, 
the property $p_h(s)\prec p_h(t)$ holds if and only if
$h(s)=\min\{h(r):s\wedge t\leq r\leq s\vee t\}$. 

\smallskip
\rem The function $h$ is not determined by $\t_h$. In particular, if
$\phi:[0,\sigma']\to [0,\sigma]$ is continuous and nondecreasing, and such that
$\phi(0)=0$ and $\phi(\sigma')=\sigma$, the tree coded by $h\circ \phi$
is isometric to the tree coded by $h$. This simple observation will
be useful later.

\begin{proposition}
\label{tree-coding}
Under the preceding assumptions, if $g$ is a nonnegative continuous function on $\t_h$
such that $g(\rho_h)=0$
and $g$ is nondecreasing with respect to the genealogical order on $\t_h$, 
the subordinate tree $(\t_h^{(g)},d_h^{(g)})$ of $(\t_h,d_h)$ with respect to the 
function $g$ is isometric to the tree $(\t_G,d_G)$ coded by the function $G=g\circ p_h$.
\end{proposition}

\proof Note that the function $G$ is nonnegative and continuous on $[0,\sigma]$,
and $G(0)=G(\sigma)=0$. We can therefore make sense of the tree $(\t_G,d_G)$ and
as above we denote the canonical projection from $[0,\sigma]$ onto $\t_G$
by $p_G$. We first notice that, for every $s,t\in[0,\sigma]$, the property 
$p_h(s)=p_h(t)$ implies $p_G(s)=p_G(t)$. Indeed, if $p_h(s)=p_h(t)$, 
then, for every $r\in [s\wedge t,s\vee t]$, we have
$p_h(s)\prec p_h(r)$ and therefore $g(p_h(s))\leq g(p_h(r))$, so that
$$G(s)=G(t)=\min_{r\in [s\wedge t,s\vee t]} G(r),$$
and $p_G(s)=p_G(t)$. We can thus write $p_G=\bp \circ p_h$, where the
function $\bp:\t_h \la \t_G$ is continuous and onto. 

Then, let $a,b\in \t_h$
and write $a=p_h(s)$ and $b=p_h(t)$, with $s,t\in[0,\sigma]$. 
We note that, for every $r\in[s\wedge t,s\vee t]$, $p_h(s)\wedge p_h(t)\prec p_h(r)$, and
furthermore $p_h(s)\wedge p_h(t)=p_h(r_0)$, if $r_0$
is any element of $[s\wedge t,s\vee t]$ at which $h$ attains its minimum over $[s\wedge t,s\vee t]$.
It follows that
$$g(p_h(s)\wedge p_h(t))=\min_{r\in[s\wedge t,s\vee t]} G(r).$$
Hence,
$$d^{(g)}_h(a,b)= g(p_h(s))+g(p_h(t))- 2 g(p_h(s)\wedge p_h(t))
= G(s)+G(t) - 2\min_{r\in[s\wedge t,s\vee t]} G(r)= d_G(s,t).$$
If  $\pi^{(g)}$ is the 
projection from $\t_h$ onto the subordinate tree $\t_h^{(g)}$, we see
in particular that  the condition $\pi^{(g)}(a)=\pi^{(g)}(b)$
implies $p_G(s)=p_G(t)$ and therefore $\bp(a)=\bp(b)$. 

It follows that $\bp= I \circ \pi^{(g)}$, where 
$I:\t_h^{(g)} \la \t_G$ is onto. It remains 
to verify that 
$I$ is isometric, but this is immediate from
the identities in the last display. \endproof

\rem It is known that any compact $\R$-tree can be represented in the form 
$\t_h$ for some function $h$ (see \cite[Corollary 1.2]{Duq}). Thus Proposition \ref{tree-coding}
provides an alternative proof of Proposition \ref{tree-quotient}.

\section{L\'evy trees}
\label{sec:Levy-tree}

In the next sections, we will consider the case where $\t$ is the (random) tree coded
by a Brownian excursion distributed under the It\^o excursion measure
$\bn(\cdot)$, and we will identify certain subordinate trees as L\'evy trees. In this section, we
recall the basic facts about L\'evy trees that will be needed later. We refer to
\cite{DLG0,DLG} for more details.

We consider a nonnegative function $\psi$ defined on $[0,\infty)$ of the type
\begin{equation}
\label{branch-meca}
\psi(r)=\alpha r + \beta r^2 + \int_{(0,\infty)} \pi(\D u)\,(e^{-ur}-1 +ur),
\end{equation}
where $\alpha,\beta\geq 0$, and $\pi(\D u)$ is a $\sigma$-finite measure on
$(0,\infty)$ such that $\int (u\wedge u^2)\,\pi(\D u) <\infty$. With any such function $\psi$,
we can associate a continuous-state branching process (see \cite{Grey} and references
therein), and $\psi$ is then called the branching mechanism function of this
process. Notice that the conditions on $\psi$ are not the most general ones, because
we restrict our attention to the critical or subcritical case. Additionally, we will
assume that 
\begin{equation}
\label{extinct}
\int_1^\infty \frac{\D r}{\psi(r)} <\infty.
\end{equation}
This condition, which implies that at least one of the two properties $\beta >0$
or $\int u\,\pi(\D u)=\infty$ holds, is equivalent to the a.s. extinction of the continuous-state branching process 
with branching mechanism $\psi$ \cite{Grey}. Special cases include $\psi(r)=r^\gamma$
for $1<\gamma\leq 2$. 

Under the preceding assumptions, one can make sense of the L\'evy tree that describes
the genealogy of the continuous-state branching process 
with branching mechanism $\psi$. We consider, under a probability measure $\mathbf{P}$, a spectrally positive L\'evy process $X=(X_t)_{t\geq 0}$
with Laplace exponent $\psi$, meaning that $\mathbf{E}[\exp(-\lambda X_t)]=\exp(t\psi(\lambda))$ for every
$t\geq 0$ and $\lambda>0$. We define the associated height process by setting, for every $t\geq 0$,
$$H_t=\lim_{\ve \to 0} \frac{1}{\ve} \int_0^t \D s\,\mathbf{1}_{\{X_s<\inf\{X_r:s\leq r\leq t\}+\ve\}}\,,$$
where the limit holds in probability under $\mathbf{P}$. Then \cite[Theorem 1.4.3]{DLG0} the process $(H_t)_{t\geq 0}$
has a continuous modification, which we consider from now on. We have $H_t=0$
if and only if $X_t=I_t$, where $I_t=\inf\{X_s:0\leq s\leq t\}$ is the past minimum process of $X$

Let $\mathbf{N}$ stand for the (infinite) excursion measure of $X-I$. Here the normalization
of $\mathbf{N}$ is fixed by saying that the local time at $0$ of $X-I$ is the process $-I$.
Let $\chi$ stand for the duration of the excursion under $\mathbf{N}$.
The height process $H$ is well defined (and has continuous paths)
under $\mathbf{N}$, and we have $H_0=H_\chi=0$, $\mathbf{N}$ a.e. 
To simplify notation, we will write $\max H=\max\{H_s:0\leq s\leq \chi\}$ under $\mathbf{N}$.

By definition (see \cite[Definition 4.1]{DLG}), the L\'evy tree with branching mechanism
$\psi$ (or in short the $\psi$-L\'evy tree) is the random compact $\R$-tree $\t_H$ coded by the function $(H_t)_{0\leq t\leq \chi}$
under $\mathbf{N}$, or more generally any random tree with the same 
distribution -- note that the distribution of the L\'evy tree is an infinite measure.
We refer to \cite{DLG0,DLG} for several results explaining in which sense the L\'evy tree
codes the genealogy of the continuous-state branching process 
with branching mechanism $\psi$. In the special case $\psi(r)=r^2/2$, $X$ is
just a standard linear Brownian motion, $H_t=2(X_t-I_t)$ is twice a reflected Brownian motion,
and the L\'evy tree is the tree coded by (twice) a positive Brownian excursion under the 
(suitably normalized)
It\^o measure. Conditioning on $\chi=1$ then yields the Brownian continuum random tree. 
When $\psi(r)=r^\gamma$ with $1<\gamma<2$, one gets the stable tree with
index $\gamma$. 

The distribution of the height of a L\'evy tree is given as follows. For every $h>0$,
$$\mathbf{N}(\mathcal{H}(\t_H)>h)=\mathbf{N}(\max H> h)= v(h),$$
where the function $v:(0,\infty)\to(0,\infty)$ is determined by
$$\int_{v(h)}^\infty \frac{\D r}{\psi(r)}= h.$$

\rem In the preceding considerations, the normalization of the infinite measure $\mathbf{N}$
is fixed by our choice of the local time at $0$ of $X-I$, and we can recover 
$\psi$ from $\mathbf{N}$ by the formulas for the distribution of $\mathcal{H}(\t_H)$ under $\mathbf{N}$. 
What happens if we multiply $\mathbf{N}$ by a constant 
$\lambda>0$\,? The tree $\t_H$ under $\lambda \mathbf{N}$ is still a L\'evy tree in the previous sense, but 
the associated branching mechanism is now $\wt\psi(r)=\lambda \psi(r/\lambda)$. 
To see this, consider the L\'evy process $X'_t=\frac{1}{\lambda}X_{\lambda t}$,
whose Laplace exponent is $\wt\psi(\lambda)$. It is not
hard to verify that the height process corresponding to $X'$ is $H'_t=H_{\lambda t}$. 
Furthermore, if $\mathbf{N}'$ is the excursion measure of $X'$ above its past minimum process,
one also checks that the distribution of $(H'_{t/\lambda})_{t\geq 0}$ under $\mathbf{N}'$
is the distribution of $(H_t)_{t\geq 0}$ under $\lambda \mathbf{N}$. However, the tree coded by
$(H'_{t/\lambda})_{t\geq 0}$ is the same as the tree coded by $(H'_{t})_{t\geq 0}$.
This shows that $\t_H$ under $\lambda \mathbf{N}$ is a L\'evy tree with branching mechanism
$\wt\psi$. Note that this is consistent with the formula for $\mathbf{N}(\mathcal{H}(\t_H)>h)$.

\medskip

We now state two results that will be important for our purposes. We first mention that, for every $h\geq 0$,
one can define under $\mathbf{P}$ a local time process of $H$ at level $h$, which is denoted by $(L^h_t)_{t\geq 0}$
and is such that the following approximation holds for every $t>0$,
\begin{equation}
\label{approx-LT}
\lim_{\ve\to 0} \mathbf{E}\Big[ \sup_{s\leq t}\Big| \frac{1}{\ve}\int_0^s \mathbf{1}_{\{h< H_s\leq h+\ve\}}\D r - L^h_s\Big|\Big] =0\,,
\end{equation}
and the latter convergence is uniform in $h$ (see \cite[Proposition 1.3.3]{DLG0}). 
When $h=0$ we have simply $L^0_t=-I_t$.
The definition of 
$(L^h_t)_{t\geq 0}$ also makes sense under $\mathbf{N}$, with a similar approximation. 

We fix $h>0$, and
let $(u_j,v_j)_{j\in J}$ be the collection of all excursion intervals of $H$ above level $h$: these are
all intervals $(u,v)$, with $0\leq u<v$, such that $H_u=H_v=h$ and 
$H_r>h$ for every $r\in (u,v)$. For each such excursion interval $(u_j,v_j)$, we define the corresponding
excursion by $H^{(j)}_s=H_{(u_j+s)\wedge v_j} -h$, for every $s\geq 0$. Then $H^{(j)}$ is a random element
of the space $C(\R_+,\R_+)$ of all continuous functions from $\R_+$ into $\R_+$. We also let
$\mathbf{N}^\circ$ be the $\sigma$-finite measure on $C(\R_+,\R_+)$ which is the ``law'' of $(H_s)_{s\geq 0}$
under $\mathbf{N}$.

\begin{proposition}
\label{excu-above}
{\rm (i)} Under the probability measure $\mathbf{P}$, the point measure
$$\sum_{j\in J} \delta_{(L^h_{u_j},H^{(j)})}(\D \ell,\D \omega)$$
is Poisson with intensity
$\mathbf{1}_{[0,\infty)}(\ell)\,\D \ell\,\mathbf{N}^\circ(\D \omega)$.

\noindent{\rm (ii)} Under the probability measure $\mathbf{N}(\cdot\mid \max H >h)$
and conditionally on $L^h_\chi$, the point measure
$$\sum_{j\in J} \delta_{(L^h_{u_j},H^{(j)})}(\D \ell,\D \omega)$$
is Poisson with intensity
$\mathbf{1}_{[0,L^h_\chi]}(\ell)\,\D \ell\,\mathbf{N}^\circ(\D \omega)$.
\end{proposition}

See \cite[Proposition 3.1, Corollary 3.2]{DLG} for a slightly more precise version of this proposition. 
It follows from the preceding proposition that L\'evy trees satisfy a branching property 
analogous to the classical branching property of Galton--Watson trees. To state this
property, we introduce some notation. If $(\t,d)$ is a (deterministic) compact $\R$-tree and $0<h< \mathcal{H}(\t)$,
we can consider the subtrees of $\t$ above level $h$. Here, a subtree above level $h$ is just the closure of a
connected component of $\{a\in\t:d(\rho,a)> h\}$. Such a subtree is itself  viewed as a rooted $\R$-tree
(the root is obviously the unique point at height $h$ in the
subtree). 

\begin{proposition}
\label{branching-Levy}
Let $\t$ be a random compact $\R$-tree defined under an infinite measure $\mathcal{N}$,
such that $\mathcal{N}(\mathcal{H}(\t)=0)=0$ and $0<\mathcal{N}(\mathcal{H}(\t)>h)<\infty$
for every $h>0$. For every $h,\ve>0$, write $M(h,h+\ve)$ for the number 
of subtrees of $\t$ above level $h$ with height greater than $\ve$. 
\begin{enumerate}
\item[\rm(i)] Suppose that $\t$ is a L\'evy tree. Then,
for every $h,\ve>0$, for every integer $p\geq 1$, the distribution 
under $\mathcal{N}(\cdot\mid M(h,h+\ve)=p)$
of the unordered collection formed by the $p$ subtrees of $\t$ above level $h$ with height greater than $\ve$
is the same as that of the unordered collection of $p$ independent 
copies of $\t$ under $\mathcal{N}(\cdot\mid \mathcal{H}(\t)>\ve)$. 
\item[\rm(ii)] Conversely, if the property stated in {\rm(i)} holds, then $\t$ is a L\'evy tree.
\end{enumerate}
\end{proposition}

The property stated in (i) is called the branching property. The fact that it holds for L\'evy trees
is a straightforward consequence of Proposition \ref{excu-above} (ii). For the converse, we refer
to \cite[Theorem 1.1]{Weill}.
Note that the branching property remains valid if we multiply
the underlying measure $\mathcal{N}$ by a positive constant, which is consistent
with the remark above.

\smallskip

To conclude this section, let us briefly comment on points of
infinite multiplicity of the L\'evy tree $\t_H$.
The multiplicity of a point $a$ of $\t_H$ is the number of connected components of 
$\t_H\backslash\{a\}$, and $a$ is called a leaf if it has multiplicity one. Suppose that there is no quadratic part in $\psi$, meaning that the constant $\beta$
in \eqref{branch-meca} is $0$ (note that condition \eqref{extinct} then implies that
$\pi$ has infinite mass). Then \cite[Theorem 4.6]{DLG}, all points of $\t_H$ have multiplicity $1$, $2$ or $\infty$. The set of all points of infinite multiplicity is a countable dense subset of $\t_H$,
and these points are in one-to-one correspondence with local minima of $H$, or with jump times of 
$X$. More precisely, 
let $a$ be a point of infinite multiplicity of $\t_H$. Then, if $s_1=\min p_H^{-1}(a)$ and $s_2=\max p_H^{-1}(a)$, we have
$s_1<s_2$,
$H_{s_1}=H_{s_2}=\min\{H_s:s_1\leq s\leq s_2\}$,
and 
$p_H^{-1}(a)=\{s\in[s_1,s_2]: H_s= H_{s_1}\}$ is a Cantor set contained in $[s_1,s_2]$. 
In terms of the L\'evy process $X$, $s_1$ is a jump time of $X$ and 
$s_2=\inf\{s\geq s_1: X_s\leq X_{s_1-}\}$, and $p_H^{-1}(a)$
consists exactly of those $s\in [s_1,s_2]$ such that
$\inf\{X_r:r\in[s_1,s]\}=X_s$. Furthermore,
if for every $r\in[0,\Delta X_{s_1}]$ we set $\eta_r=\inf\{s\geq s_1: X_{s}<X_{s_1}-r\}$,
the points of $p_H^{-1}(a)$ are all of the form $\eta_r$ or $\eta_{r-}$. The quantity $\Delta X_{s_1}$
is the ``weight'' of the point of infinite multiplicity $a$. 
See
\cite{DLG} for more details.

\section{Subordination by the Brownian snake maximum}
\label{sec:subord-maxi}

In this section, we prove Theorem \ref{subor-brown}. We start with a brief
presentation of the Brownian snake. We refer to
\cite[Chapter IV]{Zurich} for more details. 

We let
$\ww$ be the space of all finite paths in $\R$. Here a finite path is simply 
a continuous mapping $\w:[0,\zeta]\la \R$, where
$\zeta=\zeta_{(\w)}$ is a nonnegative real number called the 
lifetime of $\w$. The set $\ww$ is a Polish space when equipped with the
distance
$$d_\ww(\w,\w')=|\zeta_{(\w)}-\zeta_{(\w')}|+\sup_{t\geq 0}|\w(t\wedge
\zeta_{(\w)})-\w'(t\wedge\zeta_{(\w')})|.$$
The endpoint (or tip) of the path $\w$ is denoted by $\wh \w=\w(\zeta_{(\w)})$.
For every $x\in\R$, we set $\ww_x=\{\w\in\ww:\w(0)=x\}$. We also identify the
trivial path of $\ww_x$ with zero lifetime with the point $x$. 

The standard (one-dimensional) Brownian snake with initial point $x$ is the continuous Markov  process 
$(W_s)_{s\geq 0}$ taking values in $\ww_x$, whose distribution
is characterized by the following properties:
\begin{enumerate}
\item[(a)] The process $(\zeta_{(W_s)})_{s\geq 0}$ is a reflected Brownian motion in $\R_+$
started from $0$. To simplify notation, we write $\zeta_s=\zeta_{(W_s)}$ for every $s\geq 0$.
\item[(b)] Conditionally on $(\zeta_{s})_{s\geq 0}$, the process $(W_s)_{s\geq 0}$ is time-inhomogeneous
Markov, and its transition kernels are specified as
follows. If $0\leq s\leq
s'$,
\begin{enumerate}
\item[$\bullet$] $W_{s'}(t)=W_{s}(t)$ for every $t\leq m_\zeta(s,s'):=\min\{\zeta_r:s\leq r\leq s'\}$;
\item[$\bullet$] the random path $(W_{s'}(m_\zeta(s,s')+t)-W_{s'}(m_\zeta(s,s')))_{0\leq t\leq
\zeta_{s'}- m_\zeta(s,s')}$ is independent of $W_s$ and distributed as a
real Brownian motion started at $0$ and stopped at time 
$\zeta_{s'}- m_\zeta(s,s')$.
\end{enumerate}
\end{enumerate}

Informally, the value $W_s$ of the Brownian snake at time $s$
is a one-dimensional Brownian path started from $x$, with lifetime $\zeta_s$. 
As $s$ varies, the lifetime $\zeta_s$ evolves like reflected Brownian motion in $\R_+$. When $\zeta_s$ decreases,
the path is erased from its tip, and when $\zeta_s$ increases, the path 
is extended by adding ``little pieces'' of Brownian paths at its tip.

For every $x\in \R$, we write $\P_x$ for the probability measure under which $W_0=x$, and $\N_x$
for the (infinite) excursion measure of $W$ away from $x$. Also 
$\sigma:=\sup\{s>0: \zeta_s> 0\}$ stands for the duration of the excursion under $\N_x$. Under $\N_x$, $(\zeta_s)_{s\geq 0}$ 
is distributed according to the It\^o excursion measure $\bn(\cdot)$, and
the normalization is fixed by the formula
$$\N_x\Big(\max_{0\leq s\leq \sigma}\zeta_s >\ve\Big)= \frac{1}{2\ve}.$$

The following property of the Brownian snake will be used in several places below. Recall that $\wh W_s=W_s(\zeta_s)$ is the tip of the path $W_s$.  We say that $r\in[0,\sigma)$ is a time of right increase of $s\mapsto \zeta_s$, resp. of $s\mapsto \wh W_s$, if there
exists $\ve\in(0,\sigma-r]$ such that $\zeta_u\geq \zeta_r$, resp. $\wh W_u\geq \wh W_r$, for every $u\in[r,r+\ve]$. We can similarly define
points of left increase. Then according to \cite[Lemma 3.2]{LGP}, $\N_x$ a.e., no time $r\in(0,\sigma)$
can be simultaneously a time of (left or right) increase of $s\mapsto \zeta_s$ and a point
of (left or right) increase of $s\mapsto \wh W_s$.

Le us fix $x=0$ and argue under $\N_0$. 
As previously, we write $\t_\zeta$ for the (random) tree coded by the function $(\zeta_s)_{0\leq s\leq \sigma}$,
$p_\zeta:[0,\sigma]\la \t_\zeta$ for the canonical projection, and $\rho_\zeta=p_\zeta(0)$ for the root
of $\t_\zeta$. Properties of the Brownian snake
show that the condition $p_\zeta(s)=p_\zeta(s')$ implies $W_s=W_{s'}$, and thus we can define $W_a$ for every $a\in \t_\zeta$ by 
setting $W_a=W_s$ if $a=p_\zeta(s)$. We note that, if $a=p_\zeta(s)$ and $t\in[0,\zeta_s]$,
$W_s(t)$ coincides with $\wh W_b$ where  $b$ is the unique point of $\llbracket \rho_\zeta,a\rrbracket$  at distance $t$ from $\rho_\zeta$. For every 
$\w\in \ww$, set
$$\ov \w:=\max_{0\leq t\leq \zeta_{(\w)}} \w(t).$$
Then, the function $a\mapsto \ov W_a$
is continuous and nondecreasing on $\t_\zeta$ (if $a=p_\zeta(s)$ and $b=p_\zeta(s')$, the condition $a\prec b$
implies that $\zeta_s\leq \zeta_{s'}$ and that $W_a$ is the restriction 
of $W_b$ to the interval $[0,\zeta_s]$, so that obviously $\ov W_a\leq \ov W_b$). 
As in Theorem \ref{subor-brown}, we write $\wt\t_\zeta$ for the subordinate tree of $\t_\zeta$
with respect to the function $a\mapsto \ov W_a$.

\medskip

\noindent{\it Proof of Theorem \ref{subor-brown}.} We first verify that the 
branching property  stated in Proposition \ref{branching-Levy} holds for $\wt \t_\zeta$
under $\N_0$, and to this end we
rely on the special Markov property 
of the Brownian snake.

Let $h>0$, and, for every $\w\in \ww$, set $\tau_h(\w)=\inf\{t\in[0,\zeta_{(\w)}]:\w(t)\geq h\}$. 
Let $(a_i,b_i)_{i\in I}$ be the connected components of the open set 
$\{s\geq 0: \tau_h(W_s)<\zeta_s\}$. For every such connected component $(a_i,b_i)$,
for every $s\in[a_i,b_i]$, the path $W_s$ coincides with $W_{a_i}=W_{b_i}$
up to time $\tau_h(W_{a_i})=\zeta_{a_i}=\zeta_{b_i}=\tau_h(W_{b_i})$
(these assertions are straightforward consequences of the properties of the Brownian snake,
and we omit the details). We then
set, for every $s\geq 0$,
$$W^{(i)}_s(t) := W_{(a_i+s)\wedge b_i}(\zeta_{a_i}+t),\quad \hbox{for }
0\leq t\leq \zeta^{(i)}_s:=\zeta_{(a_i+s)\wedge b_i}-\zeta_{a_i}.$$
We view $W^{(i)}$ as a random element of the space of all
continuous functions from $\R_+$ into $\ww_h$, and the $W^{(i)}$'s are
called  the excursions 
of $W$ outside the domain $(-\infty,h)$  (see the appendix below
for further details in a more general setting). 

By a compactness argument,
only finitely many of the excursions $W^{(i)}$ hit $(h+\ve,\infty)$. Let $M_{h,\ve}$
be the number of these excursions. It follows from
Corollary \ref{SMPusual} in the appendix that, for every
$p\geq 1$, conditionally on $\{M_{h,\ve}=p\}$, the unordered
collection formed by the excursions 
of $W$ outside $(-\infty,h)$ that hit $(h+\ve,\infty)$ is distributed 
as the unordered collection of $p$ independent copies of
$W$ under $\N_h(\cdot\mid \sup\{\wh W_s:s\geq 0\} > h+\ve)$.
On the other hand, noting that $\wt\t_\zeta$ is the tree coded by
the function $[0,\sigma]\ni s\mapsto \ov W_s$ (by Proposition \ref{tree-coding}), we also see that
subtrees of $\wt\t_\zeta$ above level $h$ with height greater than $\ve$ are
in one-to-one correspondence with excursions 
of $W$ outside $(-\infty,h)$ that hit $(h+\ve,\infty)$, and if
a subtree $\wt\t^{(i)}$ corresponds to an excursion $W^{(i)}$,
$\wt\t^{(i)}$ is obtained from the excursion $W^{(i)}$ (shifted so
that it starts from $0$) by exactly the same procedure that allows
us to construct $\wt\t_\zeta$ from $W$ under $\N_0$: To be specific, $\wt\t^{(i)}$ is coded by
the function $s\mapsto \ov W^{(i)}_s-h$ just as $\wt\t_\zeta$
is coded by $s\mapsto \ov W_s$. The preceding
considerations show that $\wt\t_\zeta$ satisfies the branching property,
and therefore is a L\'evy tree.

To get the formula for $\psi_0$, we note that the
distribution of the height of $\wt \t_\zeta$ is given by
$$\N_0(\mathcal{H}(\wt \t_\zeta)>r)=\N_0\Big( \sup_{s\geq 0} \wh W_s >r\Big) =\frac{3}{2r^2},$$
where the last equality can be found in \cite[Section VI.1]{Zurich}.
Since we also know that the function $v(r)=\N_0(\mathcal{H}(\wt \t_\zeta)>r)$
solves $v'=-\psi_0(v(r))$, the formula for $\psi_0$ follows. \hfill$\square$

\section{Approximating a L\'evy tree by embedded Galton--Watson trees}

In this section, we come back to the general setting of Section \ref{sec:Levy-tree}. Our goal
is to prove that the L\'evy tree $\t_H$ is (under the probability measure
$\mathbf{N}(\cdot\mid \mathcal{H}(\t_H)>h)$ for some $h>0$)  the almost sure limit of a sequence of
embedded Galton--Watson trees, and that this limit is consistent with the
order structure of the L\'evy tree. We refer to \cite{probasur} for basic facts about Galton--Watson trees.
A key property for us is the fact that Galton--Watson trees are rooted ordered 
(discrete) trees, also called plane trees, so that there
is a lexicographical ordering on vertices.

In what follows, we argue under the probability measure $\mathbf{P}$.
Recall that $X$ is under $\mathbf{P}$ a L\'evy process with Laplace exponent $\psi$,
and that $H$ is the associated height process.
We fix an integer $n\geq 1$, and, for every integer $j\geq 0$, we consider the sequence
of all excursions of $H$ above level $j\,2^{-n}$ that hit level $(j+1)2^{-n}$.
We let 
$$0\leq \alpha^n_0<\alpha^n_1<\alpha^n_2<\cdots$$
be the ordered sequence consisting of all the initial times of these excursions, for all 
values of the integer $j\geq 0$ (so, $\alpha^n_0$ corresponds to the beginning 
of an excursion of $H$ above $0$ that hits $2^{-n}$, $\alpha^n_1$ may be either 
the beginning of an excursion of $H$ above $0$ that hits $2^{-n}$ or 
the beginning of an excursion of $H$ above $2^{-n}$ that hits $2\times2^{-n}$, and
so on). For every $j\geq 0$, we also let $\beta^n_j$ be the terminal time of the
excursion starting at time $\alpha^n_j$.

We then set, for every integer $k\geq 0$,
$$H^n_k= 2^n\, H_{\alpha^n_k}.$$

\begin{proposition}
\label{discre-height}
The process $(H^n_k)_{k\geq 0}$ is the discrete height process of a 
sequence of independent Galton--Watson trees with the same 
offspring distribution $\mu_n$.
\end{proposition}

Recall that the discrete height process of a sequence of Galton--Watson trees
gives the generation of the successive vertices in the trees, assuming that these vertices
are listed
in lexicographical order in each tree and one tree after another. See \cite{probasur} or \cite[Section 0.2]{DLG0}. The (finite) height sequence of a single tree is defined analogously.

\smallskip
\proof By construction, $\alpha^n_0$ is the initial time 
of the first excursion of $H$ above $0$ that hits $2^{-n}$. 
Notice that this excursion is distributed as $H$ under $\bN(\cdot\mid \max H \geq 2^{-n})$. Let $K\geq 1$
be the (random) integer such that $\alpha^n_{K-1}<\beta^n_0\leq \alpha^n_K$, so that
$ \alpha^n_K$ is the initial time 
of the second excursion of $H$ above $0$ that hits $2^{-n}$. 

\begin{figure}[!h]
 \begin{center}
 \includegraphics[width=14cm]{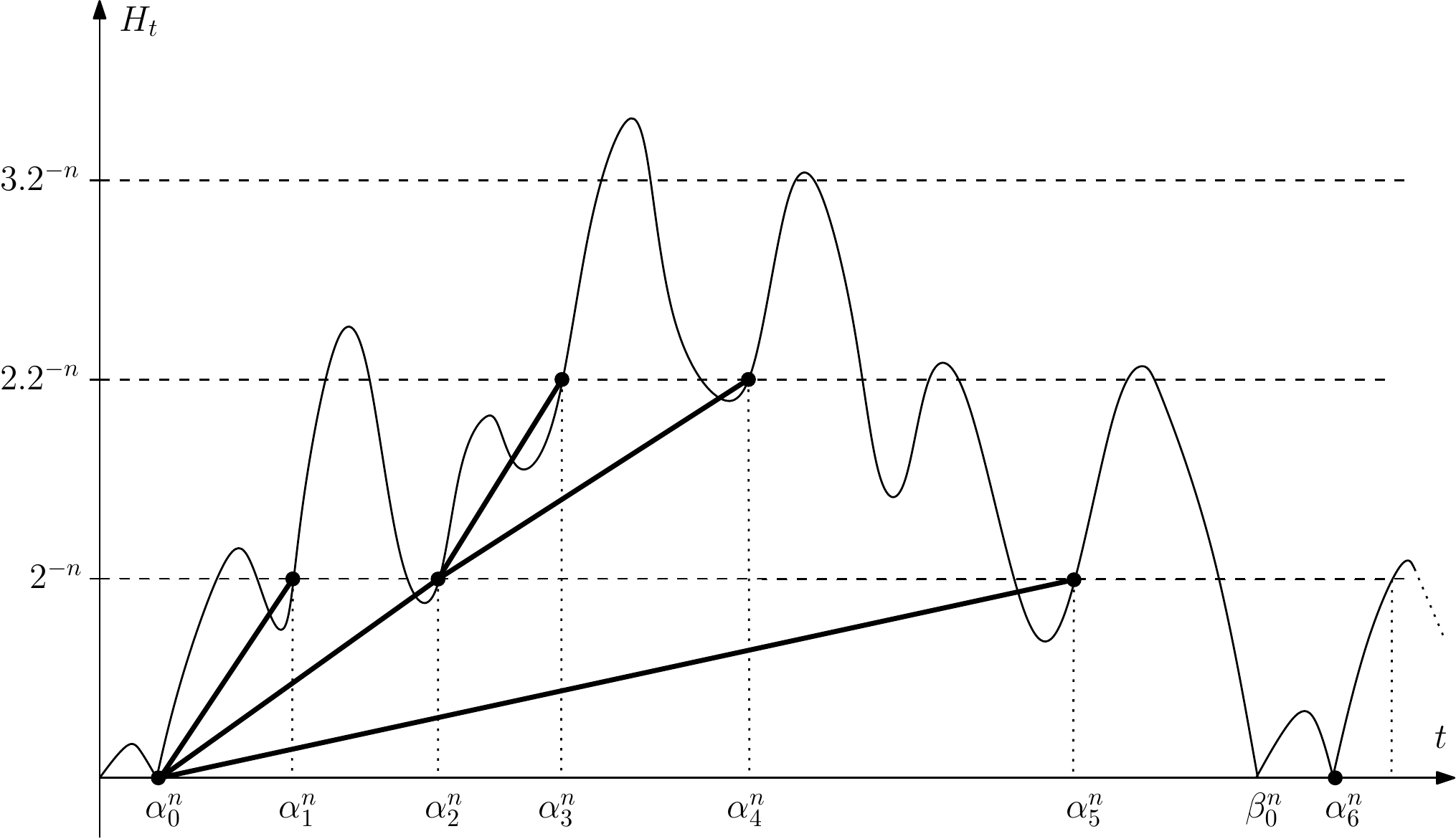}
 \caption{\label{tree-exc}
The sequence $\alpha^n_0,\alpha^n_1,\ldots$ and the tree $\t^n_0$ (in thick lines).}
 \end{center}
 \vspace{-6mm}
 \end{figure}

With the excursion of $H$ during interval $[\alpha^n_0,\beta^n_0]$, we can associate a
(plane) tree $\t^n_0$ constructed as follows. The children of the ancestor correspond to
the excursions of $H$ above level $2^{-n}$, during the time interval $[\alpha^n_0,\beta^n_0]$, that hit $2\times 2^{-n}$
and the order on these children is obviously given by the chronological order. 
Equivalently, the children of the ancestor correspond to
the indices $i\in\{1,\ldots,K-1\}$ such that $H^n_i=1$. Then, assuming 
that the ancestor has at least one child (equivalently that $K\geq 2$), the children 
of the first child of the ancestor correspond to
the excursions of $H$ above level $2\times 2^{-n}$, during the time interval $[\alpha^n_1,\beta^n_1]$, that hit $3\times 2^{-n}$,
and so on. See Fig.\ref{tree-exc} for an illustration. 

Write $N^n_0$ for the number of children of the ancestor in $\t^n_0$. It follows from Proposition \ref{excu-above} (ii) that, conditionally on $N^n_0$,
the successive excursions of $H$ above level $2^{-n}$, during the time interval $[\alpha^n_0,\beta^n_0]$, that hit $2\times 2^{-n}$
are independent and distributed as $H$ under $\bN(\cdot\mid \max H \geq 2^{-n})$ (recall that our definition 
shifts excursions above a level $h$ so that they start from $0$). Recalling the construction of the tree $\t^n_0$, we now obtain that, conditionally on $N^n_0$, the subtrees 
of $\t^n_0$ originating from the children of the ancestor are independent and distributed according to
$\t^n_0$. This just means that $\t^n_0$ is a Galton--Watson tree, and its offspring distribution $\mu_n$ is the law 
under $\bN(\cdot\mid \max H \geq 2^{-n})$ of
the number of excursions of $H$ above level $2^{-n}$ that hit $2\times 2^{-n}$. 

With the second excursion of $H$ above $0$ that hits $2^{-n}$, we can similarly associate 
a Galton--Watson tree $\t^n_1$ with offspring distribution $\mu_n$, and so on.
The trees $\t^n_0,\t^n_1,\ldots$ are independent as a consequence of the
strong Markov property of the L\'evy process $X$. By construction,
the process $(H^n_k)_{k\geq 0}$ is the discrete height process of the
sequence $\t^n_0,\t^n_1,\ldots$. 
\endproof

\begin{proposition}
\label{conv-discre}
For every $n\geq 1$, set $v_n:= 2^nv(2^{-n})=2^n \bN(\max H \geq 2^{-n})$. Then, for every $A>0$,
$$\sup_{t\leq A} | 2^{-n} H^n_{\lfloor v_n t\rfloor} - H_t| \build{\la}_{n\to\infty}^{} 0$$
in probability under $\mathbf{P}$.
\end{proposition}

\proof Recall that, for every $h\geq 0$, $(L^h_t)_{t\geq 0}$ denotes the local time of
$H$ at level $h$. It will be convenient to
introduce, for every $n\geq 1$ and every $j\geq 0$, the increasing process
$$\nn^n_{(j)}(t):=\#\{k\geq 0: H^n_k=j,L^{j2^{-n}}_{\alpha^n_k}\leq t\}.$$
As a consequence of Proposition \ref{excu-above}(i) applied with $h=j 2^{-n}$, we get that $\nn^n_{(j)}(t)$ is a
Poisson process with parameter
$v(2^{-n})=\bN(\max H \geq 2^{-n})$. 

We claim that, for every $A>0$,
\begin{equation}
\label{tech-approdis}
\lim_{n\to\infty}
\sup_{j\geq 0}\Bigg( \mathbf{E}\Big[\sup_{s\leq A} |v(2^{-n})^{-1}\,\#\{k: H^n_k=j,\alpha^n_k\leq s\} - L^{j2^{-n}}_s|\Big] \Bigg)=0.
\end{equation}

To see this, 
first observe that, for every $s>0$,
$$\nn^n_{(j)}((L_s^{j2^{-n}})-)\leq 
\#\{k: H^n_k=j,\alpha^n_k\leq s\} 
\leq\nn^n_{(j)}(L_s^{j2^{-n}}),$$
and then write
\begin{align*}
&\mathbf{E}\Big[\sup_{s\leq A} |v(2^{-n})^{-1}\,\#\{k: H^n_k=j,\alpha^n_k\leq s\} - L^{j2^{-n}}_s|\Big]\\
&\qquad\leq \sum_{p=1}^\infty \mathbf{E}\Big[\mathbf{1}_{\{p-1\leq L^{j2^{-n}}_A\leq p\}} \sup_{t\leq p} |v(2^{-n})^{-1}\,\nn^n_{(j)}(t) - t|\Big]\\
&\qquad\leq \sum_{p=1}^\infty \mathbf{P}(p-1\leq L^{j2^{-n}}_A)^{1/2}\,\mathbf{E}\Big[\sup_{t\leq p} |v(2^{-n})^{-1}\,\nn^n_{(j)}(t) - t|^2\Big]^{1/2}.
\end{align*}
Then, if $\nn(t)$ stands for a standard Poisson process, we have by a classical martingale inequality
\begin{align*}
\mathbf{E}\Big[\sup_{t\leq p} |v(2^{-n})^{-1}\,\nn^n_{(j)}(t) - t|^2\Big]&= v(2^{-n})^{-2}\mathbf{E}\Big[\sup_{t\leq v(2^{-n})p} (\nn(t)-t)^2\Big]\\
&\leq 4v(2^{-n})^{-2}\,\mathbf{E}[(\nn(v(2^{-n})p)-v(2^{-n})p)^2]\\
&=4v(2^{-n})^{-1}p.
\end{align*}
It follows that, for every $j\geq 0$,
$$\mathbf{E}\Big[\sup_{s\leq A} |v(2^{-n})^{-1}\,\#\{k: H^n_k=j,\alpha^n_k\leq s\} - L^{j2^{-n}}_s|\Big]
\leq \Big(  \sum_{p=1}^\infty (p\,\mathbf{P}(p-1\leq L^{j2^{-n}}_A))^{1/2}\Big) \times 2v(2^{-n})^{-1/2},$$
and the proof of \eqref{tech-approdis} is completed by noting that 
$v(2^{-n})\la \infty$ as $n\to\infty$, and that 
$$\sum_{p=1}^\infty (p\,\mathbf{P}(p-1\leq L^{j2^{-n}}_A))^{1/2}
\leq \sum_{p=1}^\infty (p\,\mathbf{P}(p-1\leq L^{0}_A))^{1/2}<\infty,$$
because $L^{j2^{-n}}_A$ is
bounded above in distribution by $L^0_A$ (cf Definition 1.3.1 in \cite{DLG0}), and we know that
$L^0_A=-I_A$ has exponential moments.. 

 Let $\ell\geq 1$ be an integer. By summing 
 the convergence in \eqref{tech-approdis} over possible choices
of $0\leq j< \ell 2^n$, we also obtain that
\begin{equation}
\label{tech-appro2}
\lim_{n\to\infty} \mathbf{E}\Big[\sup_{s\leq A} \Big|2^{-n}v(2^{-n})^{-1}\,\#\{k: H_{\alpha^n_k}< \ell,\alpha^n_k\leq s\} - 2^{-n}\sum_{j=0}^{\ell 2^n-1}L^{j2^{-n}}_s\Big|\Big] =0.
\end{equation}
On the other hand, we have, for every $s\geq 0$,
\begin{align*}
\Big| \int_0^s \D r\,\mathbf{1}_{\{H_r\leq \ell\}} - 2^{-n}\sum_{j=0}^{\ell 2^n-1}L^{j2^{-n}}_s\Big|
&= \Big|\sum_{j=0}^{\ell 2^n-1}\Big( \int_0^s \D r\,\mathbf{1}_{\{j2^{-n}<H_r\leq (j+1)2^{-n}\}} - 2^{-n}L^{j2^{-n}}_s\Big)\Big|\\
&\leq 2^{-n}  \sum_{j=0}^{\ell 2^n-1} \Big| 2^n \int_0^s \D r\,\mathbf{1}_{\{j2^{-n}<H_r\leq (j+1)2^{-n}\}} - L^{j2^{-n}}_s\Big|,
\end{align*}
and it follows from \eqref{approx-LT} that
\begin{equation}
\label{tech-appro3}
\lim_{n\to\infty} \mathbf{E}\Big[\sup_{s\leq A} \Big| \int_0^s \D r\,\mathbf{1}_{\{H_r\leq \ell\}} - 2^{-n}\sum_{j=0}^{\ell 2^n-1}L^{j2^{-n}}_s\Big|\Big] =0.
\end{equation}
By combining \eqref{tech-appro2} and \eqref{tech-appro3}, we get
$$
\lim_{n\to\infty} \mathbf{E}\Big[\sup_{s\leq A} \Big| \int_0^s \D r\,\mathbf{1}_{\{H_r\leq \ell\}} - v_n^{-1}\,\#\{k: H_{\alpha^n_k}< \ell,\alpha^n_k\leq s\} \Big|\Big] =0,
$$
where $v_n=2^nv(2^{-n})$. Since $\mathbf{P}(\max\{H_s:0\leq s\leq A\}\geq \ell)$ can be made arbitrarily small by choosing $\ell$
large, we have obtained that 
$$\sup_{s\leq A} \Big| v_n^{-1}\,\#\{k: \alpha^n_k\leq s\} - s\Big| \build{\la}_{n\to\infty}^{} 0$$
in probability. Elementary arguments show that this implies
\begin{equation}
\label{conv-time-change}
\sup_{t\leq A} |\alpha^n_{\lfloor v_nt\rfloor} - t|\build{\la}_{n\to\infty}^{} 0
\end{equation}
in probability, and therefore also
$$\sup_{t\leq A} |H_{\alpha^n_{\lfloor v_nt\rfloor}} - H_t|\build{\la}_{n\to\infty}^{} 0$$
in probability. This completes the proof. \endproof

\smallskip
In what follows, we will need an analog of the preceding two propositions for the height process under 
its excursion measure $\bN$. Let us fix $n\geq 1$. Under the measure $\bN(\cdot \cap\{ \max H\geq 2^{-n}\})$,
we define 
$$0= \bar\alpha^n_0<\bar\alpha^n_1<\bar\alpha^n_2<\cdots<\bar\alpha^n_{\bar{K}_n-1}$$
as the ordered sequence consisting of the initial times of 
all excursions of $H$ above level $j\,2^{-n}$ that hit level $(j+1)2^{-n}$, for all 
values of the integer $j\geq 0$.
The analog of Proposition \ref{discre-height} says that, under $\bN(\cdot \mid \max H\geq 2^{-n})$,
the finite sequence 
$$\bar H^n_k:= 2^n\, H_{\bar\alpha^n_k},\quad 0\leq k\leq \bar K_n-1$$
is distributed as the height sequence of a Galton--Watson tree with offspring distribution $\mu_n$. 
This is immediate from the fact that an excursion with distribution $\bN(\cdot \mid \{ \max H\geq 2^{-n}\})$
is obtained by taking (under $\mathbf{P}$) the first excursion of $H$ with height greater than $2^{-n}$. 
By convention, we take $\bar \alpha^n_k=\chi$ and $\bar H^n_k=0$ if $k\geq \bar K_n$.

We next fix a sequence $(n_p)_{p\geq 1}$ such that both \eqref{conv-time-change} and the convergence of Proposition \ref{conv-discre}
hold $\mathbf{P}$ a.s. along this sequence, for each $A>0$. From now on, we consider only values of
$n$ belonging to this sequence. We claim that we have then also
\begin{equation}
\label{conv-discreH}
\sup_{t\geq 0} | 2^{-n} \bar H^n_{\lfloor v_n t\rfloor} - H_t| \build{\la}_{n\to\infty}^{} 0,\qquad \bN\hbox{ a.e.}
\end{equation}
To see this, note that it suffices to argue under $\bN(\cdot\mid \max H >\delta)$ for some $\delta>0$, and
then to consider the first excursion (under $\mathbf{P}$) of $H$ away from $0$ with height greater than
$\delta$. We abuse notation by still writing $0= \bar\alpha^n_0<\bar\alpha^n_1<\bar\alpha^n_2<\cdots<\bar\alpha^n_{\bar{K}_n-1}$ for the finite sequence of times defined as explained above, 
now relative to this first excursion with height greater than $\delta$.

We observe that, provided $n$ is large enough so that $2^{-n}<\delta$, we have
$$\bar\alpha^n_k=\alpha^n_{d_n+k}, \qquad 0\leq k< \bar K_n=r_n-d_n$$
where $d_n$ is the index such that $\alpha^n_{d_n}=:d_{(\delta)}$ is the initial time of the first excursion of $H$ away from $0$ with height greater than
$\delta$, and $r_n$ is the first index $k>d_n$ such that $\alpha^n_k$
does not belong to the interval $[d_{(\delta)},r_{(\delta)}]$ associated with this excursion. Notice that $\alpha^n_{r_n}$ decreases to $r_{(\delta)}$ as $n\to\infty$. Our claim \eqref{conv-discreH} then reduces to verifying that
$$\sup_{t\geq 0} | H_{\alpha^n_{(d_n+\lfloor v_n t\rfloor)\wedge r_n}} - H_{(d_{(\delta)}+t)\wedge r_{(\delta)}}| \build{\la}_{n\to\infty}^{} 0,\qquad \mathbf{P}\hbox{ a.s.}$$
which follows from Proposition \ref{conv-discre} and \eqref{conv-time-change}, recalling that both these convergences
hold a.s. on the sequence of values of $n$ that we consider.

\section{The coding function of the subordinate tree as a time-changed height process}

In this section, we prove Theorem \ref{subordi-coding}.
We consider the Brownian snake under its excursion measure $\N_0$, and we recall the
notation $\ov W_s=\max\{W_s(t):0\leq t\leq \zeta_s\}$, and 
$\ov W_a=\ov W_s$ if $a=p_\zeta(s)$. As in Theorem \ref{subor-brown}, we write $\wt\t_\zeta$ for the subordinate tree
of the Brownian tree $\t_\zeta$ with respect to the function $a\mapsto \ov W_a$. 
\smallskip

\noindent{\it Proof of Theorem \ref{subordi-coding}.}
We set $W^*=\max\{\ov W_s:0\leq s\leq \sigma\}$. For fixed $n\geq 1$, we define a discrete plane tree
$\t^{(n)}$ under the probability measure $\N_0(\cdot\mid W^*\geq 2^{-n})$ in the following way.
The children of the root correspond to the excursions of $W$ outside $(-\infty,2^{-n})$ that
hit $2\times 2^{-n}$ (recall the definition of these excursions from the proof
of Theorem \ref{subor-brown}). Note that these excursions in turn correspond to the excursions
of the real-valued process $s\to \ov W_s$ above level $2^{-n}$ that hit $2\times 2^{-n}$
(the point is that, if $(u,v)$ is the time interval associated with an excursion of $W$ outside $(-\infty,2^{-n})$,
the process $s\mapsto \ov W_s$ remains strictly above $2^{-n}$ on the whole interval $(u,v)$). We obviously
order the children of the root according to the chronological order of the Brownian snake.

By the special Markov property, in the form given in Corollary \ref{SMPexcursion} in the appendix below, conditionally
on the number of excursions of $W$ outside $(-\infty,2^{-n})$ that
hit $2\times 2^{-n}$, these excursions listed in chronological order are
independent and distributed according to $\N_0(\cdot\mid W^*\geq 2^{-n})$, modulo
the obvious translation by $-2^{-n}$. We can thus continue the construction of the 
tree $\t^{(n)}$ by induction, and this random plane tree is a Galton--Watson tree
since it satisfies the branching property at the first generation.

Let $\mu_n$ be the offspring distribution found in Proposition \ref{discre-height} 
in the case where $\psi(r)=\psi_0(r)$. We claim that $\mu_n$ is also
the offspring distribution of $\t^{(n)}$. To see this, observe that $\mu_n$ is, by definition,
the distribution of the number of points of a $\psi_0$-L\'evy tree at height $2^{-n}$
that have descendants at height $2\times 2^{-n}$ (conditionally on the event that the height of
the tree is at least $2^{-n}$). Thanks to Theorem \ref{subor-brown} and to the fact that
$\wt\t_\zeta$ is the tree coded by the
function $s\mapsto\ov W_s$, we know that this is the same as the conditional distribution 
of the number of excursions of $s\mapsto \ov W_s$ above level $2^{-n}$ that hit $2\times 2^{-n}$, under
$\N_0(\cdot\mid W^*\geq 2^{-n})$. 

Let 
$$0\leq \xi^n_0<\xi^n_1<\xi^n_2<\cdots<\xi^n_{ K_n-1}$$
be the ordered sequence consisting of the initial times of 
all excursions of $s\to \ov W_s$ above level $j\,2^{-n}$ that hit level $(j+1)2^{-n}$, 
for all 
values of the integer $j\geq 0$. Note that each such excursion corresponds to a vertex of the tree $\t^{(n)}$, and so
$ K_n$ is just the total progeny of $\t^{(n)}$. By convention, we also define $\xi^n_{K_n}=\sigma$. Set
$$\wt H^n_k= 2^n\,\ov W_{\xi^n_k},\quad \hbox{if }0\leq k< K_n,$$
and $\wt H^n_k=0$ if $k\geq  K_n$.
Then $(\wt H^n_k,0\leq k< K_n)$ is the height sequence 
of $\t^{(n)}$ (note that the lexicographical ordering on vertices of $\t^{(n)}$ corrresponds
to the chronological order on the associated excursion initial times). Hence $(\wt H^n_k)_{k\geq 0}$
has the same distribution as the sequence $(\bar H^n_k)_{k\geq 0}$ which was defined at the end of the 
preceding section from the height process $H$ under $\bN(\cdot\mid \max H \geq 2^{-n})$.

But in fact more in true: the whole collection of the discrete sequences $(\wt H^n_k)_{k\geq 0}$
for all $n\geq 1$ has the same distribution under $\N_0$  as the similar collection
of sequences
$(\bar H^n_k)_{k\geq 0}$
constructed from the height process $H$ under $\bN$ (the reason is the fact that,
in both constructions, 
the tree at step $n$ can be obtained from the tree at step $n+1$ 
by the deterministic operation consisting in keeping only those
vertices at even generation that have at least one child, and 
viewing that set of vertices as a plane tree in the obvious manner). The convergence \eqref{conv-discreH} now allows us to set, for every $t\geq 0$,
$$\wt H_t=\lim_{n\to\infty} 2^{-n} \wt H^n_{\lfloor v_n t\rfloor}\,,$$
and the process $(\wt H_t)_{t\geq 0}$ is distributed as the height process of
the L\'evy tree with branching mechanism $\psi_0(r)=\sqrt{8/3}\,r^{3/2}$. The limit in the
preceding display holds uniformly in $t$, $\N_0$ a.e., provided we argue along the subsequence of values of $n$
introduced at the end of the preceding section. We set $\wt \chi=\sup\{s\geq 0:\wt H_s>0\}$.

We observe that the distribution of $(\wt H, (\wt H^n_k)_{n\geq 1,k\geq 0})$
under $\N_0$ is the same as that of $(H, (\bar H^n_k)_{n\geq 1,k\geq 0})$
under $\mathbf{N}$, and so we must have, for every $n\geq 1$ and $k\geq 0$,
$$\wt H^n_k=2^n\,\wt H_{\wt \alpha^n_k},$$
where $\wt\alpha^n_0<\wt \alpha^n_1<\cdots<\wt \alpha^n_{\wt K_n-1}$ are
the initial times of 
the excursions of $\wt H_s$ above level $j\,2^{-n}$ that hit level $(j+1)2^{-n}$, 
for all 
values of the integer $j\geq 0$, and $\wt\alpha^n_k=\wt\chi$ if $k\geq \wt K_n$. 
Notice that $\wt K_n=K_n$ because 
the height sequence of the tree $\t^{(n)}$  is $(\wt H^n_k)_{0\leq k\leq K_n-1}$. Also, if $n<m$
and $k\in\{0,1,\ldots,K_n-1\}$, $k'\in\{0,1,\ldots,K_m-1\}$, the property
$\tilde\alpha^n_k= \tilde\alpha^m_{k'}$ holds if and only if $\xi^n_k=\xi^n_{k'}$: Indeed these properties hold if and only if the vertex with index $k$ in the (lexicographical) ordering of $\t^{(n)}$
coincides with the vertex with index $k'$ in the ordering of $\t^{(m)}$, modulo the
identification of the vertex set of $\t^{(n)}$ as a subset of the vertex set of $\t^{(m)}$,
in the way explained above.

We need to verify that $\ov W_s$ can be written as a time change of $\wt H$. 
As a first step,
we notice that, for every $0\leq k<  K_n$,
$$2^n \ov W_{\xi^n_k}=\wt H^n_k= 2^n \wt H_{\wt\alpha^n_k}.$$
and so $\ov W_{\xi^n_k}=\wt H_{\wt\alpha^n_k}$. This suggests that the process $\Gamma$
in the statement of the theorem should be such that $\Gamma_{\xi^n_k}=\wt\alpha^n_k$,
for every $k\in\{0,1,\ldots,K_n-1\}$ and every $n$.

At this point, we observe that
\begin{equation}
\label{max-gap}
\max_{1\leq k\leq K_n} (\wt\alpha^n_k-\wt \alpha^n_{k-1}) \build{\la}_{n\to\infty}^{} 0,\qquad \N_0\ \hbox{a.e.}
\end{equation}
Indeed, if this property fails, a compactness argument gives two times $u,v\in[0,\wt\chi]$ with $u<v$ such that
$t\to \wt H_t$ is monotone nonincreasing on $[u,v]$. To see that this cannot occur,
we may replace $\wt H$ by the process $H$
constructed from a L\'evy process excursion $X$ as explained
in Section \ref{sec:Levy-tree}.  We then note that jumps of $X$ are dense
in $[0,\chi]$, and the strong Markov property shows that, for any jump time $s$ of $X$, for any $\ve>0$, we can find $s',s''\in[s,s+\ve]$,
with $s''>s'$, such that $H_{s''}>H_{s'}$ (use formula (20) in \cite{DLG0}, or
see the comments at the end of Section \ref{sec:Levy-tree}).

Let $s\in[0,\sigma)$ and, for every integer $n\geq 1$, let $k_n(s)\in\{0,1,\ldots,K_n-1\}$
be the unique integer such that $\xi^n_{k_n(s)}\leq s<\xi^n_{k_n(s)+1}$. We note that the
sequence $\xi^n_{k_n(s)}$ is monotone nondecreasing (this is obvious since
$(\xi^n_0,\ldots,\xi^n_{K_n})$ is a subset of $(\xi^{n+1}_0,\ldots,\xi^{n+1}_{K_{n+1}})$). It
follows that
the sequence $\wt\alpha^n_{k_n(s)}$ is also monotone nondecreasing:
Indeed, if $n<m$ and $k\in\{0,1,\ldots,K_n-1\}$,
$k'\in\{0,1,\ldots,K_m-1\}$ are such that $\xi^n_k\leq \xi^m_{k'}$,
we have automatically $\wt\alpha^n_k\leq\wt \alpha^m_{k'}$ since, writing
$\xi^n_k=\xi^m_{k^*}$, the fact that $k^*\leq k'$ implies that
$\wt\alpha^n_k=\wt\alpha^m_{k^*}\leq \wt\alpha^m_{k'}$. 

We can now set
$$\Gamma_s=\lim_{n\to\infty} \wt\alpha^n_{k_n(s)}.$$
Note that this limit will exist simultaneously for all $s\in[0,\sigma)$ outside a set of
$\N_0$-measure $0$. We also take $\Gamma_s=\wt\chi$ for all $s\geq \sigma$. 
Clearly $s\mapsto\Gamma_s$ is nondecreasing and, by
construction, the property $\Gamma_{\xi^n_k}=\wt\alpha^n_k$ holds 
for every $k\in\{0,1,\ldots,K_n\}$ and every $n$, and 
we have $\ov W_s=H_{\Gamma_s}$ when $s$ is of the form $\xi^n_k$. 
We also note that $s\la \Gamma_s$ is continuous as a consequence
of the property \eqref{max-gap}. To check right-continuity, observe that, if
$s\in[0,\sigma)$ is fixed, and $s'$ is such that $\xi^n_{k_n(s)}<s'<\xi^n_{k_n(s)+1}$
then, for every $m\geq n$, the property $\xi^m_{k_m(s')}\leq s'<\xi^n_{k_n(s)+1}$
forces $\wt \alpha^m_{k_m(s')}\leq \wt\alpha^n_{k_n(s)+1}$, hence
(letting $m$ tend to $\infty$) $\Gamma_s\leq \wt\alpha^n_{k_n(s)+1}$,
and use \eqref{max-gap}. Left-continuity is derived by a similar argument.

For $s\in[0,\sigma)$, set $s'=\lim\uparrow \xi^n_{k_n(s)}$ and $s''=\lim\downarrow \xi^n_{k_n(s)+1}$. Note that $s'\leq s\leq s''$. On one hand,
by passing to the limit $n\to\infty$ in the equality 
$$\ov W_{\xi^n_{k_n(s)}}=\wt H_{\alpha^n_{k_n(s)}},$$
we obtain that $\ov W_{s'}=\wt H_{\Gamma_s}$. On the other hand, $\ov W$ must be constant on the interval $[s',s'']$. To see this, 
we first
observe that $\ov W$ must be nonincreasing on $[s',s'']$
(otherwise there would be some $n$ and some 
$k\in\{0,1,\ldots,K_n-1\}$ such that $\xi^n_{k_n(s)}<\xi^n_k <\xi^n_{k_n(s)+1}$,
which is absurd). So we need 
to verify that there is no nontrivial interval $[s_1,s_2]$
such that $s\mapsto \ov W_s$ is both nonincreasing and 
nonconstant on $[s_1,s_2]$, and, to prove this,
we may replace nonincreasing by nondecreasing thanks to
the invariance of $\N_0$ under time-reversal. Argue by contradiction, and suppose that 
$s_1<s_2$ are such that that the event
where $0<s_1<s_2<\sigma$ and $s\mapsto \ov W_s$ is both nondecreasing and 
nonconstant on $[s_1,s_2]$ has  positive $\N_0$-measure. 
We can then find a stopping time $T$ such that, with positive $\N_0$-measure
on the latter event, we 
have $s_1<T<s_2$, $\ov W_T=\wh W_T$ and $\max\{W_T(s):0\leq s \leq \zeta_T\}$ is attained
only at $\zeta_T$ (take
$T=\inf\{s>s_1: \wh W_s\geq \ov W_{s_1}+\delta\}$, with $\delta>0$ small
enough). Using the strong Markov property of the Brownian snake,
we then find $r\in(T,s_2)$ such that $W_r$ is the restriction of $W_T$ to
$[0,\zeta_T-\ve]$, for some $\ve>0$, which implies $\ov W_r<\ov W_T$
and gives a contradiction with the fact that $s\mapsto \ov W_s$ is
nondecreasing on $[s_1,s_2]$. 

Finally, since $\ov W$ is constant on $[s',s'']$, we have 
$\ov W_s=\ov W_{s'}=\wt H_{\Gamma_s}$, which was
the desired result. This completes the proof of Theorem \ref{subordi-coding}. \hfill$\square$

\section{Applications to the Brownian map}
\label{sec:Br-map}

In this section, we discuss applications of the previous results to the Brownian map.
Analogously to \cite{MS}, we consider a version of the Brownian map with ``randomized volume'',
which may be constructed under the Brownian snake excursion measure $\N_0$ as follows.
Recall that $(\t_\zeta,d_\zeta)$ stands for the tree coded by $(\zeta_s)_{0\leq s\leq \sigma}$ and $p_\zeta$ is the canonical projection 
from $[0,\sigma]$ onto $\t_\zeta$. 
For 
$a,b\in \t_\zeta$, the ``lexicographical interval'' $[a,b]$ stands for the image under
$p_\zeta$ of the smallest interval $[s,t]$ ($s,t\in[0,\sigma]$) such that $p_\zeta(s)=a$ and $p_\zeta(t)=b$ (here we make the
convention that if $s>t$ the interval $[s,t]$ is equal to $[s,\sigma]\cup [0,t]$). 

For every $a\in\t_\zeta$, we set $Z_a=\wh W_s$,
where $s$ is such that $p_\zeta(s)=a$. In particular $Z_{\rho_\zeta}=0$. The 
random mapping $\t_\zeta\ni a \mapsto Z_a$
is interpreted as Brownian motion indexed by
the ``Brownian tree'' $\t_\zeta$. 

We then define a mapping $D^\circ:\t_\zeta\times \t_\zeta\la \R_+$ by setting
$$D^\circ(a,b)=Z_a+Z_b-2\max\Big\{\min_{c\in[a,b]} Z_c, \min_{c\in[b,a]} Z_c\Big\}.$$
For $a,b\in\t_\zeta$, we set
$a\approx b$ if and only if $D^\circ(a,b)=0$, or equivalently 
$$Z_a=Z_b=\max\Big\{\min_{c\in[a,b]} Z_c, \min_{c\in[b,a]} Z_c\Big\}.$$
One can verify that if $a\approx a'$ then $D^\circ(a,b)=D^\circ(a',b)$ for any $b\in\t_\zeta$
(the point is that, if $a\approx a'$ with $a\not =a'$, then 
necessarily $a$ and $a'$ are leaves of $\t_\zeta$, and the reals $t,t'\in[0,\sigma)$ such that $p_\zeta(t)=a$
and $p_\zeta(t')=a'$
are unique, which implies that $[a,b]\subset [a,a']\cup [a',b]$ and $[b,a]\subset [b,a']\cup [a',a]$). 

We also set
$$D^*(a,b)=\inf_{a_0=a,a_1,\ldots,a_p=b} \sum_{i=1}^p D^\circ(a_{i-1},a_i),$$
where the infimum is over all choices of the integer $p\geq 1$ and of the points
$a_1,\ldots,a_{p-1}$ of $\t_\zeta$.
If $a\approx a'$ then $D^*(a,b)=D^*(a',b)$ for any $b\in\t_\zeta$.
Furthermore one
can also prove that $D^*(a,b)=0$ if and only if $a\approx b$.
Since $D^*$ satisfies the triangle inequality, it follows that $\approx$
is an equivalence relation on $\t_\zeta$.
The Brownian map (with randomized volume) is the quotient space $\bm:=\t_\zeta/\!\approx$\,, which is equipped with the distance
induced by the function $D^*$ (with a slight abuse of notation, we still denote the induced distance by $D^*$).
We write $\Pi$ for the canonical projection from $\t_\zeta$ onto $\bm$, and $\bp=\Pi\circ p_\zeta$.
We also write $D^\circ(x,y)=D^\circ(a,b)$ if $x=\Pi(a)$ and $y=\Pi(b)$. 

In the usual construction of the Brownian map, one deals with the conditioned measure $\N_0(\cdot\mid \sigma=1)$
instead of $\N_0$, but otherwise the construction is exactly the same and we refer to \cite{Geo,Uniqueness}
for more details.

The Brownian map $\bm$ comes with two distinguished points. The first one $\rho=\bp(0)=\Pi(\rho_\zeta)$ corresponds to the root $\rho_\zeta$
of $\t_\zeta$. The second distinguished point 
is $\rho_*=\Pi(a_*)$, where $a_*$ is the (unique) point of $\t_\zeta$ at which $Z$ attains its minimum:
$$Z_{a_*}=\min_{a\in \t_\zeta} Z_a.$$
We will write $Z_*=Z_{a_*}$ to simplify notation.
The reason for considering $\rho_*$ comes from the fact that distances from $\rho_*$
have a simple expression. For any $a\in\t_\zeta$,
\begin{equation}
\label{dist-root}
D^*(\rho_*,\Pi(a))= Z_a - Z_*.
\end{equation}
The following ``cactus bound'' \cite[Proposition 3.1]{Geo} also plays an important role. Let $a,b\in \t_\zeta$, and let
$\gamma:[0,1]\to \bm$ be a continuous path in $\bm$ such that $\gamma(0)=\Pi(a)$
and $\gamma(1)=\Pi(b)$.
Then,
\begin{equation}
\label{cactus-bound}
\min_{0\leq t\leq 1} D^*(\rho_*,\gamma(t)) \leq \min_{c\in \llbracket a,b \rrbracket}
D^*(\rho_*,\Pi(c))= \min _{c\in \llbracket a,b \rrbracket} (Z_c-Z_*),
\end{equation}
where we recall that $\llbracket a,b \rrbracket$ is the geodesic segment between 
$a$ and $b$ in $\t_\zeta$, not to be confused with the interval $[a,b]$. 
In other words, any continuous path from $\Pi(a)$ to $\Pi(b)$ must come at least as close
to $\rho_*$ as the (image under $\Pi$ of the) geodesic segment from $a$ to $b$
in $\t_\zeta$.

We now introduce the metric net, in the terminology of \cite{MS}. For every 
$r\geq 0$, we consider the ball $B(r)$
defined by
$$B(r)=\{x\in\bm: D^*(\rho_*,x)\leq r\}.$$
For $0\leq r<D^*(\rho_*,\rho)=-Z_*$, we define the hull $B^\bullet(r)$ as
the complement of the connected component of $B(r)$ that contains $\rho$. 
Informally, $B^\bullet(r)$ is obtained from $B(r)$ by ``filling in'' the holes
of $B(r)$ except for the one containing $\rho$. Write 
$\partial B^\bullet(r)$ for the topological boundary 
of $B^\bullet(r)$. We define the {\it metric net} $\mm$ 
as the closure in $\bm$ of the union
$$\bigcup_{0\leq r<-Z_*} \partial B^\bullet(r).$$
Our goal is to investigate the structure of $\mm$. 

If $a\in\t_\zeta$ and $s\in[0,\sigma]$ are such that $p_\zeta(s)=a$,
we write $W_a=W_s$ (as previously) 
and we use the notation
$$\un W_a=\un W_s =\min\{W_s(t):0\leq t\leq \zeta_s\}=\min\{Z_b:b\in\llbracket \rho_\zeta,a\rrbracket\},$$
where the last equality holds because, as already mentioned, the quantities $W_s(t)$ for
$0\leq t\leq \zeta_s$ correspond to the values of $Z_b=\wh W_b$ along the geodesic
segment $\llbracket \rho_\zeta,a\rrbracket$. 

We then introduce the closed subset
of $\t_\zeta$ defined by
$$\Theta=\{a\in\t_\zeta: \un W_a=\wh W_a\}.$$
We note that points of $\Theta$ have multiplicity either $1$ or $2$ in $\t_\zeta$. Indeed,
there are only countably many points of multiplicity $3$, and it is not hard to see that
these points do not belong to $\Theta$.

\begin{proposition}
\label{descri-net}
Let $x\in\bm$. Then 
$x\in \mm$ if and only if  $x=\Pi(a)$ for some $a\in\Theta$.
\end{proposition}

\proof
Fix $r\in[0,-Z_*)$ and $x\in\bm$. We claim that $x\in \partial B^\bullet(r)$
if and only if we can write $x=\Pi(a)$ with both $\wh W_a=Z_*+r$, and
$$W_a(t)>Z_*+r\,,\quad\forall t\in[0,\zeta_{(W_a)}).$$
Indeed, if these conditions hold, we have $D^*(\rho_*,x)=r$ by \eqref{dist-root}, and the image 
under $\Pi$ of the geodesic segment from $a$ to $\rho_\zeta$
provides a path from $x$ to $\rho$ that stays outside $B(r)$ except
at the initial point $x$. It follows that $\Pi(a)$ belongs to $\partial B^\bullet(r)$.

Conversely, if $x\in \partial B^\bullet(r)$, then it is obvious that 
$D^*(\rho_*,x)=r$ giving $\wh W_a=Z_*+r$ for any $a$ such that $\Pi(a)=x$.
Write
$x=\lim x_n$, where $x_n \notin B^\bullet(r)$, and, for every $n$, let
$a_n\in \t_\zeta$ such that $\Pi(a_n)=x_n$. The fact that 
$x_n \notin B^\bullet(r)$ implies that, for every $c$ belonging to the geodesic segment between
$a_n$ and $\rho_\zeta$, we have $Z_c> Z_*+r$ (otherwise the cactus bound \eqref{cactus-bound}
would imply that any path between $x_n$ and $\rho$ visits $B(r)$, which is a contradiction). 
By compactness, we may assume that $a_n$ converges to $a$ as $n\to\infty$, and we have
$\Pi(a)=x$. We then get that the property $Z_c> Z_*+r$
holds for $c$ belonging to the geodesic segment between
$a$ and $\rho_\zeta$, except possibly for $c=a$. This completes the
proof of our claim.

It follows from the claim that the property $x=\Pi(a)$ for some $a\in\t_\zeta$ such that 
$\wh W_a = \un W_a$ holds for every $x\in\mm$ (this property holds if 
$x\in \partial B^\bullet(r)$ for some $0\leq r<-Z_*$ and is preserved under 
passage to the limit, using the compactness of $\t_\zeta$). Conversely, suppose that
this property holds, with $a\not =\rho_\zeta$
to discard a trivial case. If the path $W_a$ hits its minimum only at its
terminal point, the first part of the proof shows that $x\in \partial B^\bullet(r)$
for $r= \wh W_a-Z_*$. If the path $W_a$ hits its minimum
both at its terminal time and at another time, then Lemma 16 in \cite{ALG}
shows that $a=\lim a_n$, where $\wh W_{a_n}<\wh W_a$ for every $n$. Then the image
under $\Pi$ of the
first point $b_n$ on the ancestral line of $a_n$ such that $\wh W_{b_n}= \wh W_{a_n}$
belongs to $\partial B^\bullet(r_n)$, with $r_n= \wh W_{a_n}-Z_*$. Noting that
$b_n$ must lie on the geodesic segment between $a$ and $a_n$
in the tree $\t_\zeta$, we see that we have
also $a=\lim b_n$, so that we get that $x=\Pi(a)=\lim \Pi(b_n)$ belongs to $\mm$. 
\endproof

\rem The preceding arguments are closely related to \cite[Section 3]{Hull} (see in particular 
formula (16) in \cite{Hull}), which deals with the slightly different setting of the Brownian plane.

\smallskip

If $x\in \mm$ and $a\in \Theta$ is such that $x=\Pi(a)$, the image under $\Pi$ of the geodesic segment from 
$a$ to $\rho_\zeta$  provides a path in $\bm$ that stays in the complement of $B(r)$ for every
$0\leq r< D^*(\rho^*,x)=Z_a-Z_*$
(indeed the values of $Z_b$ for $b$ belonging to this segment are the numbers
$W_a(t)\geq \un W_a=\wh W_a=Z_a$). It follows that all points belonging to a geodesic from $x$
to $\rho_*$ also belong to $\mm$.

We note that we can define
an ``intrinsic'' metric on $\mm$ by setting, for every $x,y\in\mm$,
\begin{equation}
\label{deltastar}
\Delta^*(x,y)=\build{\inf_{x=x_0,x_1,\ldots,x_k=y}}_{x_1,\ldots,x_{k-1}\in\mm}^{}
\sum_{i=1}^k D^\circ(x_{i-1},x_i).
\end{equation}
It is obvious that $\Delta^*(x,y)\geq D^*(x,y)$. In particular,
$\Delta^*(x,y)=0$ implies $x=y$, and it follows that 
$\Delta^*$ is  a metric on $\mm$. The quantity 
$\Delta^*(x,y)$ corresponds to the infimum of the lengths
(computed with respect to $D^*$) of paths from $x$ to $y$ that
are obtained by the concatenation of pieces of geodesics from points
of $\mm$ to $\rho^*$ (we already noticed that these geodesics stay in $\mm$).
We have clearly $\Delta^*(\rho_*,x)=D^*(\rho_*,x)=D^\circ(\rho_*,x)$ for every $x\in \mm$,
and $\Delta^*$-geodesics from $x$ to $\rho_*$ 
coincide with $D^*$-geodesics from $x$ to $\rho_*$ (if $x,y\in\mm$ and $x,y$ belong to the
same $D^*$-geodesic to $\rho_*$, the results of \cite{Geo} imply
that $D^*(x,y)=D^\circ(x,y)=\Delta^*(x,y)$).

\smallskip

\noindent{\bf Remark.} 
The topology induced by $\Delta^*$
on $\mm$ coincides with the topology induced by $D^*$. Since 
$\Delta^*\geq D^*$ and $(\mm,D^*)$
is compact, it is enough to prove that 
$(\mm,\Delta^*)$
is also compact. However, if $(x_n)_{n\geq 1}$ is a sequence in $\mm$,
we may write $x_n=\Pi(a_n)$, with $a_n\in\Theta$, and then extract a 
subsequence $(a_{n_k})$ that converges to $a_\infty$ in $\t_\zeta$. We have
$a_\infty\in \Theta$ because $\Theta$ is closed. Furthermore the fact that
$a_{n_k}$ converges to $a_\infty$ implies that $D^\circ(a_{n_k},a_\infty)$
tends to $0$, and therefore $\Delta^*(x_{n_k},\Pi(a_\infty))$ also tends to $0$,
showing that $(x_n)_{n\geq 1}$
has a convergent subsequence in $(\mm,\Delta^*)$. 

\smallskip

The preceding proposition shows that the metric net $\mm$ has close connections
with the subset $\Theta$ of $\t_\zeta$.
The latter set is itself related to the subordinate tree of $\t_\zeta$ 
with respect to the function $a\mapsto -\un W_a$. By Theorem \ref{subordi-coding}
(and an obvious symmetry argument)
we can construct a
process $(H_t)_{0\leq t\leq \chi}$ distributed as the height process of the L\'evy tree
with branching mechanism $\psi_0(r)=\sqrt{8/3}\,r^{3/2}$, and a
continuous random process $(\Gamma_s)_{s\geq 0}$ 
with nondecreasing sample paths such that $\Gamma_0=0$, $\Gamma_\sigma=\chi$,
and for every $s\in[0,\sigma]$,
\begin{equation}
\label{repremin}
-\un W_s= H_{\Gamma_s}. 
\end{equation}

We define a random equivalence $\sim$ relation on $[0,\chi]$, by requiring that the graph of 
$\sim$ is the smallest closed symmetric subset of $[0,\chi]^2$ that contains 
all pairs $(s,t)$ with $s\leq t$, $H_s=H_t$, and $H_r>H_s$ for
all $r\in(s,t)$. We leave it to the reader to check that this set is indeed the graph
of an equivalence relation (use the comments at the end of Section \ref{sec:Levy-tree}). In addition to the pairs $(s,t)$ satisfying the previous relation, the 
graph of $\sim$ contains a countable collection of pairs $(u,v)$, each of them
associated with a point of infinite multiplicity $a$ of the tree $\t_H$
by the relations $u=\min p_H^{-1}(a)$ and $v=\max p_H^{-1}(a)$.

We denote the quotient space $[0,\chi]/\!\sim$ by $\ll$. Then $\ll$ can be identified with the  ``looptree'' associated
with $H$. Roughly speaking (see \cite{CK} for more details) the looptree is
obtained by replacing each point of infinite multiplicity $a$ of the tree $\t_H$
by a loop of ``length'' equal to the weight of $a$, so that the subtrees that are the
connected components of the complement of $a$ in the tree branch along this loop in
an order determined by the coding function $H$. Note that the looptree associated
with $H$ is equipped in \cite{CK} with a particular metric. Here we avoid introducing this
metric on $\ll$, because it will be more relevant to our applications to
introduce a pseudo-metric that will be described below.

Let us introduce the right-continuous inverse of $\Gamma$. For every $u\in [0,\chi)$,
we set
$$\tau_u:=\inf\{s\geq 0: \Gamma_s>u\}.$$
By convention, we also set $\tau_\chi=\sigma$. The left limit
$\tau_{u-}$ is equal to $\inf\{s\geq 0:\Gamma_s=u\}$, and $\Gamma$
is constant on every interval $[\tau_{u-},\tau_u]$. Note that
$\Gamma_{\tau_u}=u$ and thus $-\un W_{\tau_u}=H_u$. 

We next consider the subset $\Theta^1$ of $\Theta$ defined as follows. 
If
$a\in \Theta$, we say that $a\in \Theta^1$ 
if there exist $s\in[0,\sigma)$ and $\ve\in(0,\sigma-s)$ such that $p_\zeta(s)=a$
and 
the function $r\mapsto \un W_r$ is constant on the interval $[s,s+\ve]$.
Notice that only leaves (points $a$ of multiplicity $1$, for which there 
is a single value of $s$ with $p_\zeta(s)=a$) may belong to $\Theta^1$.
Indeed, if $a\in\Theta$ has multiplicity two, and $s_1,s_2$ are the two elements 
of $[0,\sigma)$ such that $p_\zeta(s_1)=p_\zeta(s_2)=a$, then 
both $s_1$ and $s_2$ are times of (left or right) increase of $\zeta$, and,
together with the property $\wh W_a=\un W_a$, this implies that 
$\wh W_{s}$ takes values strictly less than $\un W_a$ 
immediately after $s_1$, resp. immediately after $s_2$
(see Section \ref{sec:subord-maxi}).
 
 One can describe the elements of
$\Theta^1$ in the following way. Let $b\in \Theta$
such that $b$ has a strict descendant $c$ with $\un W_c=\un W_b$ (in the
terminology of \cite{ALG}, $b$ is an excursion debut above the minimum). Let
$[s_1,s_2]\subset [0,\sigma]$ be the interval whose image under
$p_\zeta$ gives all descendants of $b$. Then the set
$$\{r\in[s_1,s_2]: \un W_r =\un W_{s_1}\hbox{ and }\wh W_r>\un W_{s_1}\}$$
is an open subset of $[s_1,s_2]$, and the
(image under $p_\zeta$ of the) left end of each of its connected components
belongs to $\Theta^1$. Furthermore any element of
$\Theta^1$ can be obtained in this way.

We set $\Theta^\circ=\Theta\backslash \Theta^1$.
\begin{lemma}
\label{tech-rep}
For $u,v\in[0,\chi]$, the property $u\sim v$ holds if and only if $p_\zeta(\tau_u)=p_\zeta(\tau_v)$.
Furthermore, the mapping $\Psi: [0,\chi]\ni u\mapsto p_\zeta(\tau_u)$ induces a bijection 
from $\ll=[0,\chi]/\!\sim$ onto $\Theta^\circ$, which will be denoted by $\Phi$. 
\end{lemma}

\proof We first show that $p_\zeta(\tau_u)$ belongs to $\Theta^\circ$,
for every $u\in[0,\chi)$. By the definition of
$\tau_u$, $\Gamma_{\tau_u+\ve}>\Gamma_{\tau_u}$ for every $\ve>0$,
and this implies that $\un W$ is not constant on $[\tau_u,\tau_u+\ve]$
(use \eqref{repremin} and the fact that $H$ is not constant on any nontrivial interval). In
particular, we must have $\un W_{\tau_u}=\wh W_{\tau_u}$ and
therefore $p_\zeta(\tau_u)\in \Theta$. The fact that $p_\zeta(\tau_u)\in \Theta^\circ$
is then immediate from the definition of $\Theta^1=\Theta\backslash\Theta^\circ$, recalling that points of $\Theta^1$ are
leaves of $\t_\zeta$.

We then verify that, if $a\in \Theta^\circ$, there exists $u\in[0,\chi)$ with $p_\zeta(\tau_u)=a$.
We can write $a=p_\zeta(s)$ where the mapping $r\mapsto \un W_r$ is not constant on
$[s,s+\ve]$ for every $\ve >0$. Using the formula $-\un W_r= H_{\Gamma_r}$, it follows
that $\Gamma_{s+\ve}>\Gamma_s$ for every $\ve >0$, and thus $s=\tau_{\Gamma_s}$.
Finally $a=p_\zeta(\tau_u)$ with $u=\Gamma_s$.

Next let us prove that $u\sim v$ implies $p_\zeta(\tau_u)=p_\zeta(\tau_v)$. Let $u\sim v$ and without loss of generality suppose that $0<u<v<\chi$. 
We first assume
that $H_u=H_v$ and $H_r>H_u$ for every $r\in(u,v)$. Then, we must
have $\un W_s<\un W_{\tau_u}=\un W_{\tau_{v-}}$ for every $s\in (\tau_u,\tau_{v-})$
-- note that $\un W$ is constant over any interval $[\tau_{r-},\tau_r]$. This implies
that $\zeta_s\geq \zeta_{\tau_u}=\zeta_{\tau_{v-}}$ for every $s\in (\tau_u,\tau_{v-})$ 
(if there exists $s\in (\tau_u,\tau_{v-})$ such that $\zeta_s< \zeta_{\tau_u}$, then, for every
$\delta>0$ small enough, the properties of the
Brownian snake allow us to find such an $s$ with the additional property that 
$W_s$ is the restriction of $W_{\tau_u}$ to $[0,\zeta_s-\delta]$, which contradicts $\un W_s<\un W_{\tau_u}$ -- and we can make a symmetric argument if
there exists $s\in(\tau_u,\tau_{v-})$ such that $\zeta_s< \zeta_{\tau_{v-}}$). 
It follows that $p_\zeta(\tau_u)=p_\zeta(\tau_{v-})$. Moreover, since 
$\tau_{v-}$ is a point of left increase of $\zeta$, $\tau_{v-}$ cannot be a point of right increase of
$\wh W$, so that there are values of $s> \tau_{v-}$ arbitrarily close to $\tau_{v-}$ such that
$\un W_s<\un W_{\tau_{v-}}$, and therefore $\Gamma_s>\Gamma_{\tau_{v-}}$. It follows that we have $\tau_{v-}=\tau_v$ giving $p_\zeta(\tau_u)=p_\zeta(\tau_v)$ as desired. 

Suppose then that $u\sim v$ but the property $H_r>H_u$ for every $r\in(u,v)$ does not hold.
Then $(u,v)$ is the limit of a sequence $(u_n,v_n)$ such that, for every $n$, 
$H_{u_n}=H_{v_n}$ and $H_r>H_{u_n}$ for every $r\in(u_n,v_n)$. We must have $u_n<u$ and $v_n>v$. By the first part of the argument, 
we have $\zeta_s\geq \zeta_{\tau_{u_n}}=\zeta_{\tau_{v_n}}$ for
$s\in[\tau_{u_n},\tau_{v_n}]$, and, letting $n$ tend to $\infty$, we get 
$\zeta_s\geq \zeta_{\tau_{u-}}=\zeta_{\tau_{v}}$, for every $s\in[\tau_{u-},\tau_v]$.
Then the fact that $\tau_{u-}$ is a point of right increase of $\zeta$ implies 
that $\tau_{u-}=\tau_u$ (by the same argument as above)  and we conclude again
that $p_\zeta(\tau_u)=p_\zeta(\tau_v)$. 

Finally, it remains to prove that the property $p_\zeta(\tau_u)=p_\zeta(\tau_v)$ implies $u\sim v$. Note that
$\tau_u=\tau_v$ is only possible if $u=v$, so that we may assume that $\tau_u<\tau_v$. Then $a=p_\zeta(\tau_u)=p_\zeta(\tau_v)$
is a point of multiplicity $2$ of $\t_\zeta$
(since $a\in\Theta$, $a$ cannot have multiplicity $3$ in $\t_\zeta$), and the points $p_\zeta(s)$ for $\tau_u\leq s\leq \tau_v$
are descendants of $a$, so that $\un W_s\leq \un W_a$ for every $\tau_u\leq s\leq \tau_v$. It follows that
$H_r\geq H_u$ for $u\leq r\leq v$ (write $H_r=-\un W_{\tau_r}$). If $H_r> H_u$ for $u< r< v$, this means that $u\sim v$
and we are done. Otherwise there exists $r\in(u,v)$ such that $H_r=H_u$, and this means that $a$ has a strict descendant $b=p_\zeta(\tau_r)$
such that $\un W_b=\un W_a$. This implies that the path $W_a$ hits its minimal value only at its terminal time
(otherwise $W_b$ would have two equal local minima). We know that, just before $\tau_u$, there are values of $s$
such that $\zeta_s<\zeta_{\tau_u}$ (otherwise $\tau_u$ would a time of local minimum
of $\zeta$, but this is excluded since such times 
correspond to points of multiplicity $3$ of $\t_\zeta$ and thus never satisfy $\un W_s=\wh W_s$), and it follows that there
are times $s<\tau_u$ arbitrarily close to $\tau_u$ such that $\un W_s> \un W_{\tau_u}$, and
thus $H_{\Gamma_s}<H_u$. Hence, if we set
$u_n=\sup\{r<u: H_r=H_u-\frac{1}{n}\}$, we have $u_n\la u$ as $n\to\infty$. Similarly, if
$v_n=\inf\{r>v: H_r=H_v-\frac{1}{n}\}$ we have $v_n\la v$ as $n\to\infty$. Clearly $u_n\sim v_n$
so that we also get $u\sim v$. \endproof

We write  $\mathrm{p}_\mathcal{L}$ for the
canonical projection from $[0,\chi]$ onto $\mathcal{L}=[0,\chi]/\!\!\sim$. If $\alpha\in\mathcal{L}$
and $\alpha=\mathrm{p}_\mathcal{L}(s)$, we will also write $H_\alpha=H_s$. 

In a way similar to the definition of intervals in $\t_\zeta$, we can define
intervals  in $\mathcal{L}$. If $\alpha,\beta\in \ll$,
we set $[\alpha,\beta]=\mathrm{p}_\ll([s,t])$, where $s,t\in[0,\chi]$ are such that $\mathrm{p}_\ll(s)=\alpha$
and $\mathrm{p}_\ll(t)=\beta$ and $[s,t]$ is as small as possible (here again we use the
convention $[s,t]=[s,\chi]\cup[0,t]$ if $t<s$). 

We will identify the metric net $(\mm,\Delta^*)$
with a quotient space of the looptree $\ll$. Informally, the latter
quotient space is obtained by identifying two points $\alpha$ and $\beta$
if they face each other at the same height in $\ll$: This means that
we require that $H_\alpha=H_\beta$, and that vertices ``between''
$\alpha$ and $\beta$ have a smaller height. To make this more
precise, we define, for every $\alpha,\beta\in \ll$,
$$\dd^\circ(\alpha,\beta)= 2\min\Bigg(
\max_{\gamma\in[\alpha,\beta]} H_\gamma, 
\max_{\gamma\in[\beta,\alpha]} H_\gamma\Bigg) - H_\alpha -H_\beta,$$
and 
$$\dd^*(\alpha,\beta)=
\inf_{\alpha_0=\alpha,\alpha_1,\ldots,\alpha_{k-1},\alpha_k=\beta} \;
\sum_{i=1}^k \dd^\circ(\alpha_{i-1},\alpha_i),$$
where the infimum is over all possible choices of the integer $k\geq 1$
and of $\alpha_1,\ldots,\alpha_{k-1}\in \ll$. 

The following statement is a reformulation, in a more precise form, of Theorem \ref{identi-net}
stated in the introduction.

\begin{theorem}
\label{identi-metricnet}
For $\alpha,\beta\in \ll$, set $\alpha\simeq \beta$ if and only if 
$\dd^*(\alpha,\beta)=0$. Then the property $\alpha\simeq\beta$
holds if and only if $\dd^\circ(\alpha,\beta)=0$, or equivalently
\begin{equation}
\label{identi-pairs}
H_\alpha=H_\beta=\min\Bigg(
\max_{\gamma\in[\alpha,\beta]} H_\gamma, 
\max_{\gamma\in[\beta,\alpha]} H_\gamma\Bigg).
\end{equation}
Furthermore, $\dd^*$ induces a metric on the quotient space $\ll/\!\simeq$. 
If $\Phi:\ll \la \Theta^\circ$ denotes the bijection of Lemma \ref{tech-rep},
$\Pi\circ \Phi$ induces an isometry from $(\ll/\!\simeq, \dd^*)$
onto $(\mm, \Delta^*)$. 
\end{theorem}

\noindent{\bf Remark.} It is not a priori obvious that \eqref{identi-pairs}
defines an equivalence relation on $\ll$. This property follows from the fact that
\eqref{identi-pairs} holds if and only if $\dd^*(\alpha,\beta)=0$, which we
derive in
 the following proof from the relations between $\ll$ and the Brownian map.

\proof
We first verify that, if $\alpha,\beta\in \ll$
and $a=\Phi(\alpha),b=\Phi(\beta)$, we have
\begin{equation}
\label{identi-tech1}
\dd^\circ(\alpha,\beta)= D^\circ(a,b).
\end{equation}
Let $s\in[0,\chi]$ be such that $\alpha=\mathrm{p}_\ll(s)$. Note that
we have then $a=p_\zeta(\tau_s)$ by the definition of $\Phi$. Hence,
$$H_\alpha=H_s=-\un W_{\tau_s}=-Z_{\tau_s}=-Z_a.$$
Similarly, we have $H_\beta=-Z_b$.

So the proof of \eqref{identi-tech1} reduces to checking that
$$\max_{\gamma\in[\alpha,\beta]} H_\gamma = \max_{c\in[a,b]} (-Z_c).$$
To get this equality, we write
\begin{align*}
\max_{\gamma\in[\alpha,\beta]} H_\gamma&=\min_{\mathrm{p}_\ll(s)=\alpha,\mathrm{p}_\ll(t)=\beta} \Big( \max_{r\in[s,t]} H_r\Big)\\
&= \min_{\mathrm{p}_\ll(s)=\alpha,\mathrm{p}_\ll(t)=\beta} \Big( \max_{r\in[\tau_s,\tau_t]} (-Z_r)\Big)\\
&= \min_{p_\zeta(u)=a,p_\zeta(v)=b} \Big( \max_{r\in[u,v]} (-Z_r)\Big)\\
&=\max_{c\in[a,b]} (-Z_c).
\end{align*}
The second equality holds because $H_r=-\un W_{\tau_r}=-Z_{\tau_r}$, and  $\un W$ stays constant on intervals
$[\tau_{r-},\tau_r]$. To justify the third equality, we note that the elements $u$ of $[0,\sigma]$
such that $p_\zeta(u)=a$ are exactly the reals $u=\tau_s$ where $s$ is such that $\mathrm{p}_\ll(s)=\alpha$: Since
$a\in \Theta^\circ$, any $u\in[0,\sigma]$ such that $p_\zeta(u)=a$ must be of the form $u=\tau_s$
with $s\in[0,\chi]$, and if $u$ is of this form, the property $p_\zeta(u)=a$ is equivalent to $\mathrm{p}_\ll(s)=\alpha$
by Lemma \ref{tech-rep}. This completes the proof of \eqref{identi-tech1}.

We then claim that, for every $x\in\mm$, there exists $a\in\Theta^\circ$ such that $\Pi(a)=x$. 
Let $x\in\mm$. 
By Proposition \ref{descri-net}, we know that $x=\Pi(a)$ where $a\in\Theta$.
It may happen that $a\in\Theta\backslash \Theta^\circ$, but then
we can write $a=p_\zeta(s)$ where $\un W_r=\un W_s$ for every
$r\in[s,s+\ve]$, for some $\ve>0$. The latter property if only
possible if $W_s$ hits its minimal value both at its terminal time $\zeta_s$
and at another time $\eta\in(0,\zeta_s)$. If $s'=\inf\{r>s:\zeta_r\leq \eta\}$, 
$W_{s'}$ is the restriction of $W_s$ to $[0,\eta]$.
By the results recalled in Section \ref{sec:subord-maxi}, the fact that $s'$ is a time of left increase for $\zeta$ implies that there
are values of $r<s'$ arbirarily close to $s'$ such that $\un W_r<\un W_{s'}=\un W_s$. 
If we set $\wt s=\inf\{r>s:\un W_r<\un W_s\}$, we have
$s<\wt s<s'$ and $p_\zeta(\wt s)\in\Theta^\circ$. Furthermore,
$\un W_r=\un W_s=\un W_{\tilde s}$ for every $r\in[s,\wt s]$, and it follows from the 
definition of the equivalence relation $\approx$ that $p_\zeta(\wt s)\approx p_\zeta(s)$, hence
$x=\Pi(p_\zeta(s))=\Pi(p_\zeta(\wt s))$. Our claim is proved.

Let $x,y\in \mm$, and $a,b\in \Theta^\circ$ such that $\Pi(a)=x$ and $\Pi(b)=y$.
From \eqref{deltastar} and the preceding claim, we may write 
$$\Delta^*(x,y)=\build{\inf_{a=a_0,a_1,\ldots,a_k=b}}_{a_1,\ldots,a_{k-1}\in\Theta^\circ}^{}
\sum_{i=1}^k D^\circ(a_{i-1},a_i),$$
We then use
the bijection $\Phi$ of Lemma \ref{tech-rep} to observe that, if $\alpha=\Phi^{-1}(a)$
and $\beta=\Phi^{-1}(b)$, we have also, thanks to \eqref{identi-tech1},
\begin{equation}
\label{identi-tech2}
\Delta^*(x,y)=\build{\inf_{\alpha=\alpha_0,\alpha_1,\ldots,\alpha_k=\beta}}_{\alpha_1,\ldots,\alpha_{k-1}\in\ll}^{}
\sum_{i=1}^k \dd^\circ(\alpha_{i-1},\alpha_i)= \dd^*(\alpha,\beta).
\end{equation}
In particular, if $\alpha,\beta\in \ll$ and $x=\Pi\circ \Phi(\alpha)$, $y=\Pi\circ\Phi(\beta)$, we see that the
condition $\dd^*(\alpha,\beta)=0$ holds if and only if $\Delta^*(x,y)=0$, and (since $\Delta^*(x,y)\geq D^*(x,y)$)
the latter condition holds if and only if $D^\circ(x,y)=0$, or equivalently $\dd^\circ(\alpha,\beta)=0$ (by \eqref{identi-tech1}). This gives the
first assertion of the theorem.

Then $\dd^*$ is symmetric and satisfies the triangle inequality, hence induces a metric on the quotient space $\ll/\!\simeq$. 
From the property $\Delta^*(\Pi\circ\Phi(\alpha),\Pi\circ\Phi(\beta))=  \dd^*(\alpha,\beta)$, we see that 
the relation $\alpha\simeq \beta$ implies $\Pi\circ\Phi(\alpha)=\Pi\circ\Phi(\beta)$, so that 
$\Pi\circ \Phi$ induces a mapping from $\ll/\!\simeq$ to $\mm$. This mapping is onto since 
$\Pi(\Delta^\circ)=\mm$, and is an isometry by \eqref{identi-tech2}. \endproof

\section{The holes in the metric net}

In this section, we continue our applications to the Brownian map. 
We keep the notation and assumptions of the preceding section. In particular, the
process $(H_s)_{0\leq s\leq \chi}$, which is distributed under $\N_0$ as the height process of the 
$\psi_0$-L\'evy tree, was introduced so that the representation formula \eqref{repremin} holds.

Our goal is to discuss the connected components 
of the complement of the metric net in the Brownian map (these are called Brownian disks in 
\cite{MS}). We again argue under the excursion measure $\N_0$. For every $s>0$, we denote by
$\yy_s$ the (total mass of the) exit measure of the Brownian snake from $(-s,\infty)$, see \cite[Chapter V]{Zurich}.
Then $(\yy_s)_{s>0}$ has a c\`adl\`ag modification that we consider from now on (see the discussion 
in Section 2.5 of \cite{ALG}). For every $s>0$, we can also consider the total local time 
of $H$ at level $s$, which we denote by $L^s_\chi$
in agreement with Section \ref{sec:Levy-tree}. The Ray-Knight theorem of \cite[Theorem 1.4.1]{DLG0}
shows that $(L^s_\chi)_{s>0}$ is distributed under $\N_0$ according to the excursion measure 
of the continuous-state branching process with branching mechanism $\psi_0$, and therefore
has also a c\`adl\`ag modification.

\begin{lemma}
\label{identi-exit}
We have $L^s_\chi=\yy_s$ for every $s>0$, $\N_0$ a.e.
\end{lemma}

\proof Let $s>0$ and, for every $\ve>0$, let $N_{s,\ve}$ be the number of excursions of the Brownian snake
outside $(-s,\infty)$ that hit $-s-\ve$. Then an easy application of the special Markov property gives
$$\yy_s=\lim_{\ve\to 0} v(\ve)^{-1}\,N_{s,\ve}\;,\qquad \N_0\ \hbox{a.e.}$$
where $v(\ve)=\N_0(\max \wh W_s >\ve)=\N_0(\mathcal{H}(\t_H)>\ve)$.
On the other hand, \eqref{repremin} shows that $N_{s,\ve}$ is also the 
number of excursions of $H$ above $(s,\infty)$ that hit $s+\ve$. By comparing the 
preceding approximation of $\yy_s$ with \cite[Theorem 4.2]{DLG}, we arrive at the
stated result. \endproof

\smallskip
Following closely \cite{ALG}, we say that $a\in \Theta$ is an excursion debut if
$a$ has a strict descendant $b$ such that $Z_c>Z_a$ for every $c\in \rrbracket a,b\rrbracket$.
We also say that $m\in\R_+$ is a local minimum of $H$ if there exist
$s\in(0,\chi)$ and $\ve\in(0,s\wedge(1-s))$ such that $H_s=m$ and 
$H_r\geq m$ for every $r\in (s-\ve,s+\ve)$. 

We now claim that the following sets are in one-to-one correspondence:
\begin{enumerate}
\item[(a)] The set of all connected components of $\bm\backslash \mm$.
\item[(b)] The set of all connected components of $\t_\zeta\backslash\Theta$.
\item[(c)] The set of all excursion debuts.
\item[(d)] The set of all jump times of the exit measure process $\yy$.
\item[(e)] The set of all points of infinite multiplicity of $\t_H$.
\item[(f)] The set of all local minima of $H$.
\end{enumerate}

Let us explain these correspondences. 
First the fact that local minima of $H$ correspond to 
points of infinite multiplicity of $\t_H$ was explained at the end of Section \ref{sec:Levy-tree}.
Recall that, for every point of infinite multiplicity of $\t_H$,
there is a Cantor set of local minimum times corresponding to the associated
local minimum (see the end of Section \ref{sec:Levy-tree}).
Then, by \cite[Theorem 4.7]{DLG}, each branching point $b$ 
(necessarily of infinite multiplicity) of $\t_H$ 
corresponds to a discontinuity of $s\mapsto L^s_\chi=\yy_s$ at time $H_b$,
and the corresponding jump $\Delta\yy_s$ is the weight of the branching point $b$. By 
\cite[Proposition 36]{ALG}, discontinuity times for $\yy_s$ are in one-to-one 
correspondence with excursion debuts, and a discontinuity time $s$ corresponds to an
excursion debut $a$ such that $s=-Z_a$. The fact that excursion debuts are in one-to-one correspondence
with connected components of $\t_\zeta\backslash\Theta$ is Proposition 20 in \cite{ALG}:
If $a$ is an excursion debut, the associated connected component $\cc$ is the collection
of all strict descendants $b$ of $a$ such that $Z_c>Z_a$ for every $c\in \rrbracket a,b\rrbracket$,
and the boundary $\partial \cc$ consists of all
descendants $b$ of $a$ such that $Z_b=Z_a$ and $Z_c>Z_a$ for every $c\in \rrbracket a,b\llbracket$. Furthermore, the ``boundary size'' of $\cc$ may be defined as the quantity $\Delta\yy_s$,
if $s$ is the associated jump time of the process $\yy$ (this is also the weight of the
corresponding point of infinite multiplicity of $\t_H$). 
Finally, the fact that the sets (a) and (b) are also in one-to-one correspondence 
is a consequence of the following lemma.

\begin{lemma}
\label{connec-compo}
Let $\cc$ be a connected component of $\t_\zeta\backslash\Theta$.
Then $\Pi(\cc)$ is a connected component of $\bm\backslash \mm$, and $\partial \Pi(\cc)=
\Pi(\partial \cc)$. 
\end{lemma}

\proof 
Let $a\in \t_\zeta$ be the excursion debut such that $\cc$ is 
the collection
of all strict descendants $b$ of $a$ such that $Z_c>Z_a$ for every $c\in \rrbracket a,b\rrbracket$.
We first observe that $\Pi(\cc)$ is an open subset of $\bm\backslash \mm$.
This follows from the fact that the topology of $\bm$ is the quotient topology
and $\Pi^{-1}(\Pi(\cc))=\cc$ (to derive the latter equality, note that, if
$b\in\cc$ and $b'\in \t_\zeta$ are such that $b\approx b'$, we
have $ \min\{Z_c:c\in \llbracket b,b'\rrbracket\}=Z_b=Z_{b'} > Z_a$, and it follows that
$b'\in\cc$). Since $\Pi(\cc)$ is connected, in order to get
the statement of the lemma, we need only verify that, if 
$x\in \bm \backslash \mm$ is such that there is a continuous path $(\gamma(t))_{0\leq t\leq 1}$
that stays in $\bm\backslash \mm$ and
connects $x$ to a point $y$ of $\Pi(\cc)$, then $x\in\Pi(\cc)$. 
We argue by contradiction and assume that $x\notin \Pi(\cc)$. We then
set $t_0:=\inf\{t\in(0,1]:\gamma(t)\in \Pi(\cc)\}$. Clearly, $\gamma(t_0)$
belongs to the boundary $\partial\Pi(\cc)$ of $\Pi(\cc)$. 
On the other hand, it is easy to verify that $\partial\Pi(\cc)\subset \Pi(\partial\cc)$
(if $z\in \partial\Pi(\cc)$, we can write $z=\lim \Pi(a_n)$ where $a_n\in \cc$, and,
by extracting a subsequence, we can assume that $a_n \la a_\infty$ in $\t_\zeta$,
so that we have $z=\Pi(a_\infty)$, and $a_\infty \in \partial\cc$ since 
$a_\infty\in \cc$ would imply $z\in \Pi(\cc)$, contradicting $z\in \partial\Pi(\cc)$). 
So $\gamma(t_0)\in \Pi(\partial\cc)\subset \Pi(\Theta)= \mm 
$, which contradicts our assumption that $\gamma$ stays in $\bm\backslash \mm$.

For the last assertion, it remains to see that $ \Pi(\partial \cc)\subset \partial \Pi(\cc)$. 
This is straightforward:
If $b\in \partial \cc$, we have automatically
$Z_b=Z_a$ and $b\in\Theta$, so that $\Pi(b)\in \mm$, forcing
$\Pi(b)\in \partial \Pi(\cc)$. 
 \endproof

It is worth giving a direct interpretation of the correspondence between sets (c) and (f)
above. If $a$ is an excursion debut, we can write $a=p_\zeta(s_1)=p_\zeta(s_2)$,
where $s_1<s_2$, and the image under $p_\zeta$ of the interval $[s_1,s_2]$
corresponds to the exploration of descendants of $a$ in $\t_\zeta$. It follows that
we have $\un W_c\leq \un W_a=Z_a$ for every $c\in p_\zeta([s_1,s_2])$. There are
points $b\in p_\zeta([s_1,s_2])$ such that $Z_b<Z_a$ (in fact one can find such points 
arbitrarily close to $a$, as a consequence of the fact that points of increase for $\zeta$
cannot be points of increase for $\wh W$). If $b$ is such a point we can consider the
last ancestor $c$ of $b$ such that $Z_c\leq Z_a$, noting that the definition of an excursion
debut implies that $c$ is a strict descendant of $a$. Then if $t_1<t_2$ are the two times
in $(s_1,s_2)$
such that $p_\zeta(t_1)=p_\zeta(t_2)=c$, one verifies that both $\Gamma_{t_1}$ and $\Gamma_{t_2}$
are local minimum times of $H$ corresponding to the local minimum
$-Z_a$. In fact the set of all these local minimum times consists of all $\Gamma_r$
for $r\in(s_1,s_2)$ such that $p_\zeta(r)$ belongs to the boundary of the
connected component of $\t_\zeta\backslash\Theta$ associated with $a$. 

\smallskip
We will now establish that the boundary of any connected component 
of $\bm\backslash \mm$ is a simple loop. To this end, it is 
convenient to introduce the L\'evy process excursion $(X_t)_{0\leq t\leq \chi}$
associated with $H$ (see Section \ref{sec:Levy-tree}). Note that
$X$ can be reconstructed as a measurable function of $H$, and that any
branching point of $H$ corresponds to a unique jump time of $X$, such that
the size of the jump is the weight of the branching point.

Recall the notation $\Psi$ for the mapping introduced in Lemma \ref{tech-rep}.

\begin{proposition}
\label{bdry-compo}
Let $C$ be a connected component of $\bm\backslash \mm$, 
let $a$ be the associated excursion debut, and let
$s\in(0,\chi)$ be the associated jump time of $X$. For every
$r\in[0,\Delta X_s]$, set
$$\eta_r=\inf\{t\geq s: X_t<X_s-r\}.$$
Then the mapping
$$r\mapsto \gamma(r)=\Pi\circ \Psi(\eta_r), \quad 0\leq r\leq \Delta X_s,$$
defines a simple loop in $\bm$, whose initial and end points are equal
to $\Pi(a)$. Furthermore, the range of $\gamma$ is the boundary of $C$.
\end{proposition}

\proof Write $\cc$ for the connected component of 
$\t_\zeta\backslash\Theta$ such that $\Pi(\cc)=C$.
Recall that $\partial C=\Pi(\partial\cc)$ by Lemma \ref{connec-compo}.
For every $x\in \partial C$, the fact that $x=\Pi(b)$
for some $b\in\partial \cc$ forces $Z_x=Z_b=Z_a$. 
Let $[s_1,s_2]$ be the interval corresponding to the descendants of $a$
in the coding of $\t_\zeta$. Since $s_1$ and $s_2$ are both (left or right) increase
times of $\zeta$, they cannot be increase times for $\wh W$, and this implies $s_1=\tau_{\Gamma_{s_1}}$
and $s_2=\tau_{\Gamma_{s_2}}$.
We next observe that any point of $\partial\cc$ is of the form
$p_\zeta(\tau_r)$ or $p_\zeta(\tau_{r-})$ for some $r\in[0,\chi]$
such that $s_1\leq \tau_r\leq s_2$. We just noticed that this is true for 
$a$. If $b\in \partial\cc$ and $b\not= a$, we
can write $b=p_\zeta(u)$ with $u\in(s_1,s_2)$, and Lemma  16 in \cite{ALG}
implies that there are values of $v$ arbitrarily close to $u$ such that
$\wh W_v<\wh W_u$, which implies $u=\tau_{\Gamma_u}$
or $u=\tau_{\Gamma_u-}$. 

Furthermore, we
have $\Pi(p_\zeta(\tau_{r-}))=\Pi(p_\zeta(\tau_r))$ for any $r\in(0,\chi]$:
If $\tau_{r-}<\tau_r$, the fact that $\un W$ stays constant
on the time interval $[\tau_{r-},\tau_r]$, together with the properties 
$p_\zeta(\tau_r)\in\Theta$ and $p_\zeta(\tau_{r-})\in\Theta$ (which follow from 
Lemma \ref{tech-rep}),
 implies that 
$p_\zeta(\tau_{r-})\approx p_\zeta(\tau_{r})$.
Hence, any point of $\partial C$ is of the form
$\Pi(p_\zeta(\tau_r))=\Pi\circ\Psi(r)$
for some $r\in[0,\chi]$
such that $s_1\leq \tau_r\leq s_2$, and we have then
$H_r=-\un W_{\tau_r}=-Z_{p_\zeta(\tau_r)}=-Z_a$. 

Next note that the condition $s_1\leq \tau_r\leq s_2$ 
holds if and only if
$\Gamma_{s_1}\leq r\leq \Gamma_{s_2}$, and that 
$[\Gamma_{s_1},\Gamma_{s_2}]$ is the interval corresponding
to the descendants of the branching point associated with $C$,
in the coding of the tree $\t_H$. 
Using the comments of the end of Section \ref{sec:Levy-tree}, we see that this interval
is the same as $[s,\eta_{\Delta X_s}]$ with the notation of the proposition.

Set $h=-Z_a$ to simplify notation. The preceding discussion
shows that any point of $\partial C$ is of the form
$\Pi(p_\zeta(\tau_r))$, with $r\in [s,\eta_{\Delta X_s}]$
and $H_r=h$. Conversely, for any $r$ satisfying these
conditions, we have $\Pi(p_\zeta(\tau_r))\in\partial C$: This 
follows from the fact that $p_\zeta(\tau_r)\in \partial\cc$ under these
conditions. Indeed, we have then, recalling that $p_\zeta(\tau_r)\in\Theta$,
$$Z_{p_\zeta(\tau_r)}=\un W_{\tau_r}=-H_r=Z_a,$$
and the fact that $p_\zeta(\tau_r)$ is a descendant of $a$, which belongs to
$\Theta$ and satisfies $Z_{p_\zeta(\tau_r)}=Z_a$, implies that $p_\zeta(\tau_r)\in \partial\cc$.

Note that the mapping $r\mapsto \eta_r$ is right-continuous. From
the results recalled at the end of Section \ref{sec:Levy-tree},
the times $r\in [s,\eta_{\Delta X_s}]$ such that $H_r=h$
are exactly all reals of the form $\eta_u$ or $\eta_{u-}$
for some $u\in [0,\Delta X_s]$. Moreover, for all $u$
such that $\eta_{u-}<\eta_{u}$, we have $H_r> h$
for every $r\in(\eta_{u-},\eta_u)$, so that $\eta_{u-}\sim
\eta_u$ and (by Lemma \ref{tech-rep}) $p_\zeta(\tau_{\eta_{u-}})=
p_\zeta(\tau_{\eta_u})$.

We have thus obtained that the range of the mapping
$\gamma$ of the proposition coincides with $\partial C$.
We note that $\gamma(0)=\gamma(\Delta X_s)=\Pi(a)$,
by the fact that $s=\eta_0 \sim \eta_{\Delta X_s}$ (because $s=\min p_H^{-1}(\alpha)$
and $\eta_{\Delta X_s}=\max p_H^{-1}(\alpha)$, if $\alpha=p_H(s)$ is the point
of infinite multiplicity of $\t_H$ associated with $C$). 
It remains to verify that $\gamma$ is continuous and that
the restriction of $\gamma$ to $[0,\Delta X_s)$ is one-to-one. The latter
property is easy because the values of $\eta_u$ for
$0\leq u<\Delta X_s$ correspond to distinct equivalence classes
for $\sim$, and we can use Lemma \ref{tech-rep}. Next we note 
that $\gamma$ is right-continuous with left limits by construction, and that
the left limit of $\gamma$ at $u\in(0,\Delta X_s]$ is 
$$\gamma(u-)=\Pi(p_\zeta(\tau_{(\eta_{u-})-})).$$
We already noticed that
$\Pi(p_\zeta(\tau_{r-}))=\Pi(p_\zeta(\tau_r))$ for any $r\in(0,\chi]$, 
and thus we have, for any $u\in(0,\Delta X_s]$,
$$\Pi(p_\zeta(\tau_{(\eta_{u-})-}))=\Pi(p_\zeta(\tau_{\eta_{u-}}))=\Pi(p_\zeta(\tau_{\eta_{u}}))$$
where the second equality holds because 
$p_\zeta(\tau_{\eta_{u-}})=
p_\zeta(\tau_{\eta_u})$ as mentioned above. This shows that $\gamma(u-)=\gamma(u)$
and completes the proof. \endproof

Let us conclude this section with some comments. 
Recalling that the Brownian map $\bm$ is homeomorphic to 
the two-dimensional sphere \cite{LGP}, we get
from Jordan's theorem that all connected components
of $\bm\backslash \mm$ are homeomorphic to the disk.
In fact these connected components are called {\it Brownian
disks} in \cite{MS}. If $C$ is a given connected component
of $\bm\backslash \mm$, the structure of $C$ -- in a sense
that we do not make precise here -- is described by 
the associated component $\cc$ of $\t_\zeta\backslash \Theta$,
and the values of $Z$ on $\cc$ (shifted so that the boundary values
vanish). The preceding data correspond to what is called
an excursion of the Brownian snake above its minimum 
in \cite{ALG}. One key result of \cite{ALG} states that conditionally
on the exit measure process $(\yy_s)_{s>0}$, the excursions above
the minimum are independent, and the distribution 
of the excursion corresponding to a jump $\Delta \yy_s$
is given by a certain  ``excursion measure'' conditioned
on the boundary size being equal to $\Delta\yy_s$. 
This suggests that one can reconstruct the Brownian map
by first considering the metric net $\mm$ (which is a
measurable function of $H$) and then glueing independently
on each ``hole'' of the metric net (associated with
a point of infinite multiplicity of the tree $\t_H$) a Brownian disk
corresponding to a Brownian snake excursion 
whose boundary size is the weight of the point of infinite multiplicity.
We postpone a more precise version of the previous
discussion to forthcoming work (see also \cite{MS}). 

\section{Subordination by the local time}

In this section, which is mostly independent of the previous ones, we generalize the subordination by the maximum
discussed in Section \ref{sec:subord-maxi}. To this end, we
deal with the Brownian snake associated with 
a more general spatial motion.
Specifically, we consider a strong Markov process $\xi$ with
continuous sample paths with values in $\R_+$, 
and we write $P_x$ for a probability measure under 
which $\xi$ starts from $x$. We assume that $0$ is a
regular recurrent point for $\xi$, and that
\begin{equation}
\label{time-at-0}
E_0\Big[\int_0^\infty \D t\,\mathbf{1}_{\{\xi_t=0\}}\Big]=0.
\end{equation}
We can then define the
local time process $(L(t),t\geq 0)$ of $\xi$ at $0$ (up to a multiplicative constant). We make the following
continuity assumption: there exist two reals $p>0$ 
and $\ve>0$, 
and a constant $C$ such that, for every $t\in [0,1]$ and $x\in\R_+$,
\begin{align}
\label{contxi}
E_x\Big[\Big(\sup_{r\leq t} |\xi_r-x|\Big)^p\Big]&\leq C\,t^{2+\ve},\\
\label{contlocal}
E_0[L(t)^p]&\leq C\,t^{2+\ve}.
\end{align}
We write $\mathscr{N}$ for the excursion measure of $\xi$ away from $0$
associated with the local time process $L(\cdot)$, and $\eta$
for the duration of the excursion under $\mathscr{N}$.

Under the preceding assumptions, the Brownian snake
whose spatial motion is the pair $(\xi,L)$ is defined by a straightforward
adaptation of properties (a) and (b) stated at the beginning of
Section \ref{sec:subord-maxi} (see \cite[Chapter IV]{Zurich}
for more details), and we denote this
process by $(W_s,\Lambda_s)$, where for every $s\geq 0$,
$$W_s=(W_s(t))_{0\leq t\leq \zeta_s}\ ,\quad \Lambda_s=(\Lambda_s(t))_{0\leq t\leq \zeta_s}.$$
For every $(x,r)\in \R_+\times \R_+$, let $\N_{(x,r)}$ denote the excursion measure
of $(W,\Lambda)$ away from $(x,r)$. Under $\N_{(x,r)}$, the ``lifetime process''
$(\zeta_s)_{s\geq 0}$ is distributed according to the It\^o measure $\bn(\cdot)$, and as
above we let
$\sigma:=\sup\{s\geq 0:\zeta_s>0\}$ stand for the duration of the excursion 
$(\zeta_s)_{s\geq 0}$. As previously, $\t_\zeta$ denotes the tree coded by
$(\zeta_s)_{0\leq s\leq \sigma}$
and $p_\zeta: [0,\sigma]\la \t_\zeta$ is
the canonical projection. 

We write $Y_0$ for the total mass of the exit measure of $(W,\Lambda)$ from 
$(0,\infty)\times \R_+$ (see \cite[Chapter V]{Zurich} or the appendix below
for the definition of exit measures). This makes sense under the excursion measures
$\N_{(x,r)}$ for $x>0$. 

Let $\wh \Lambda_s=\Lambda_s(\zeta_s)$ be the total local time at $0$
accumulated by the path $W_s$. If $a=p_\zeta(s)$, we write $\wh\Lambda_a=
\wh \Lambda_s$ (this does not depend on the choice of $s$). Then the 
function $a\mapsto \wh\Lambda_a$ is nondecreasing with respect to the
genealogical order. 

\begin{theorem}
\label{Levy-subor}
Under $\N_{(0,0)}$, the subordinate tree $\wt \t$ of $\t_\zeta$ with respect
to the function $a\mapsto \wh\Lambda_a$ is a L\'evy tree whose 
branching mechanism $\psi$ can be described as follows:
$$\psi(r)= 2\int m(\D x)\,u_r(x)^2,$$
where $m$ is the invariant measure of $\xi$ defined by
$$\int m(\D x)\,\varphi(x)= \mathscr{N}\Big(\int_0^\eta \D t\,\varphi(\xi_t)\Big),$$
and the function $(u_r(x))_{r\geq 0,x>0}$ is given by
$$u_r(x)= \N_{(x,0)}(1-\exp(-r Y_0)).$$
\end{theorem}

\proof As in the proof of Theorem \ref{subor-brown}, we make use of the special Markov property
of the Brownian snake. Fix $r>0$, and consider
the domain $D_r=\R_+\times [0,r)$. Write 
$\z^{D_r}$ for the exit measure from $D_r$. The first-moment
formula for exit measures \cite[Proposition V.3]{Zurich} shows that $\z^{D_r}$ is 
$\N_{(0,0)}$ a.e. supported on $\{(0,r)\}$, so that
we can write
$$\z^{D_r} = Y_r\,\delta_{(0,r)},$$
where $Y_r$ is a nonnegative random variable.
Let $\e^{D_r}$ stand for the $\sigma$-field generated by the
paths $(W_s,\Lambda_s)$ before they exit $D_r$. By Corollary
\ref{SMPusual} in the appendix, under $\N_{(0,0)}$ and conditionally on $\e^{D_r}$, the
excursions of the Brownian snake $(W,\Lambda)$ ``outside'' $D_r$
form a Poisson measure with intensity $Y_r\,\N_{(0,r)}$. 
Now notice that, for every $h>0$, subtrees of $\wt \t$ above 
level $r$ that hit ${r+h}$ correspond to those among these excursions that
exit $D_{r+h}$ (we again use Proposition \ref{tree-coding} to obtain that
$\wt \t$ is the tree coded by $s\to \wh \Lambda_s$). 
As in the proof of Theorem \ref{subor-brown} (we omit a few details here), it follows that the distribution of $\wt \t$ 
under $\N_{(0,0)}$ satisfies the branching property of Proposition
\ref{branching-Levy}, and so $\wt \t$ 
under $\N_{(0,0)}$ must be a L\'evy tree.

To determine the branching mechanism of this L\'evy tree, we fix $R>0$, and,
for $0\leq r<R$, we set
$$U_\lambda(x,r)=\N_{(x,r)}(1-\exp(-\lambda Y_R)).$$
By \cite[Theorem V.4]{Zurich}, $U_\lambda$ satisfies the integral equation
$$U_\lambda(x,r) + 2\,E_{(x,r)}\Big[ \int_0^{\tau_R} U_\lambda(\xi_s,L_s)^2 \,\D s\Big]=\lambda$$
where the Markov process $(\xi,L)$ starts from $(x,r)$ under the
probability measure $P_{(x,r)}$, and $\tau_R=\inf\{t\geq 0: L(t)\geq R\}$ 
is the exit time from $D_R$ for the process $(\xi,L)$.  When $x=0$,
excursion theory for $\xi$ gives
$$E_{(0,r)}\Big[ \int_0^{\tau_R} U_\lambda(\xi_s,L_s)^2 \,\D s\Big]= 
\int_r^R \mathscr{N}\Big(\int_0^\eta U_\lambda(\xi_t,\ell)^2\D t\Big) \D \ell
= \int_r^R \D \ell \int m(\D y)\,U_\lambda(y,\ell)^2.$$
Set $v_\lambda(x,r)=U_\lambda(x,R-r)$ for $0< r\leq R$. By a
translation argument, $v_\lambda(x,r)$ does not depend on our choice 
of $R$ provided that $R\geq r$. It follows from the preceding
considerations that
$$v_\lambda(0,r) +2 \int_0^r \D \ell \int m(\D y)\,v_\lambda(y,\ell)^2 =\lambda.$$
On the other hand, by applying the special Markov property (Corollary \ref{SMPusual})
to the domain $(0,\infty)\times \R_+$, we have, for every $y> 0$ and $r>0$,
$$v_\lambda(y,r)= \N_{(y,0)}(1-\exp(-\lambda Y_r))= \N_{(y,0)}(1-\exp(-Y_0 v_\lambda(0,r))
=u_{v_\lambda(0,r)}(y),$$
with the notation introduced in the theorem. We conclude that
\begin{equation}
\label{branmeca}
v_\lambda(0,r)+ \int_0^r \D \ell\,\psi(v_\lambda(0,\ell)) =\lambda,
\end{equation}
where $\psi$ is as in the statement of the theorem. Note that 
the functions $r\mapsto u_r(x)$ are monotone increasing, and so is $\psi$. Then \eqref{branmeca} also implies that
$v_\lambda(0,r)$ is a continuous nonincreasing function of $r$, that tends
to $\lambda$ as $r\to 0$. It follows that $\psi(r)<\infty$ for every $r>0$, and 
then by dominated convergence that $\psi$ is continuous on $[0,\infty)$. 
The unique solution of \eqref{branmeca} is given by
$$\int_{v_\lambda(0,r)}^\lambda \frac{\D \ell}{\psi(\ell)} = r
$$
(in particular, we must have $\int_{0+}\psi(\ell)^{-1}\D\ell=\infty$). As $\lambda\to\infty$,
$v_\lambda(0,r)$ converges to $\N_{(0,0)}(Y_r\not =0)$, which coincides with $\N_{(0,0)}(\mathcal{H}(\wt\t)>r)$,
where $\mathcal{H}(\wt \t)$ denotes the height of $\wt \t$. Hence, the function 
$v(r)= \N_{(0,0)}(\mathcal{H}(\wt\t)>r)$ is given by
$$\int_{v(r)}^\infty \frac{\D \ell}{\psi(\ell)} = r,$$
and this suffices to establish that the branching mechanism of $\wt\t$ is $\psi$. \endproof

The formula for $\psi$ that appears in Theorem \ref{Levy-subor} is not explicit and 
in general does not allow the calculation of this function. We will now argue
that we can identify $\psi$, up to a multiplicative constant, if 
$\xi$ satisfies a scaling property. From now on until the end of the section, 
we assume
(in addition to the previous hypotheses) that there exists a constant $\alpha>0$ such that, for every
$x\geq 0$ and  $\lambda>0$, the law of
$$(\lambda^{\alpha}\xi_{t/\lambda})_{t\geq 0}$$
under $P_x$ coincides with the law of $(\xi_t)_{t\geq 0}$ under $P_{x\lambda^{\alpha}}$.
In other words, the process $\xi$ is a self-similar Markov process with values in $[0,\infty)$, see
the survey \cite{PR} for more information on this class of processes.
A particular case  (with $\alpha=1/2$) is the Bessel process of dimension $d\in(0,2)$.

The excursion measure $\mathscr{N}$ must then satisfy a similar
scaling invariance property. More precisely, for every $\lambda>0$, the law of 
$$(\lambda^{\alpha}\xi_{t/\lambda})_{t\geq 0}$$
under $\mathscr{N}$ must be equal to $\lambda^{\beta}$ times the law of 
$(\xi_t)_{t\geq 0}$ under $\mathscr{N}$, for some constant $\beta\in(0,1)$. The fact 
that $\beta<1$ is clear since the scaling property implies that
$\mathscr{N}(\eta >r)=r^{-\beta}\mathscr{N}(\eta>1)$ and we must have
$\mathscr{N}(\eta \wedge 1)<\infty$. The inverse local time 
of $\xi$ at $0$ is then a stable subordinator of index $\beta$, which is consistent
with assumption \eqref{contlocal}.

\begin{proposition}
\label{branchmecastable}
Under the preceding assumptions, there exists a constant $c>0$
such that $\psi(r)=c\,r^{1+\beta}$. 
\end{proposition}

\proof We first observe that
$$m(\D x)=c'\,x^{-1+\frac{1-\beta}{\alpha}}\,\D x,$$
for some positive constant $c'$. To see this, we write, for every $\lambda>0$,
\begin{align*}
\int m(\D x)\,\varphi(\lambda x)&= \mathscr{N}\Big( \int_0^\eta \D t\,\varphi(\lambda \xi_t)\Big)\\
&= \lambda^{-1/\alpha}\mathscr{N}\Big( \int_0^{\lambda^{1/\alpha}\eta} \D t\,\varphi(\lambda \xi_{t/\lambda^{1/\alpha}})\Big)\\
&= \lambda^{\beta/\alpha-1/\alpha}\mathscr{N}\Big( \int_0^\eta \D t\,\varphi(\xi_t)\Big)\\
&=\lambda^{\frac{\beta-1}{\alpha}} \int m(\D x)\,\varphi(x).
\end{align*}
It follows that $m$ has the form stated above. 
We then observe that, for every $\lambda>0$, we can also consider the
following scaling transformation of the Brownian snake:
$$W'_s(t)=\lambda\,W_{\lambda^{-2/\alpha}s}(\lambda^{-1/\alpha}t),\quad \hbox{for } 0\leq t\leq \zeta'_s:=\lambda^{1/\alpha}\zeta_{\lambda^{-2/\alpha}s},$$
and the ``law'' of $W'$ under $\N_{(x,0)}$ coincides with $\lambda^{1/\alpha}$
times the ``law'' of $W$ under $\N_{(\lambda x,0)}$. Furthermore
the exit measure $Y'_0$ associated with $W'$ is equal to $\lambda^{1/\alpha} Y_0$ (we leave the details as an exercise for the reader). 
With the notation of Theorem \ref{Levy-subor}, it follows that, for every $x>0$ and $\lambda>0$,
$$u_r(\lambda x)=\N_{(\lambda x,0)}(1-\exp(-r Y_0))=\lambda^{-1/\alpha} \N_{(x,0)}(1-\exp(-r \lambda^{1/\alpha} Y_0))=
\lambda^{-1/\alpha}\,u_{\lambda^{1/\alpha}  r}(x).$$
Hence, for every $r>0$ and $\mu>0$,
$$
\psi(\mu r)= c\int \D x\,x^{-1+\frac{1-\beta}{\alpha}}\,u_{\mu r}(x)^2
=c\int \D x\,x^{-1+\frac{1-\beta}{\alpha}}\,\mu^2\,u_r(\mu^\alpha x)^2
=\mu^{1+\beta}\,\psi(r),
$$
using the change of variables $y=\mu^\alpha x$. This completes the proof. \endproof

\rem If $\xi=|B|$ is the absolute value of a linear Brownian motion $B$, then a famous
theorem of L\'evy asserts that the pair $(\xi,L)$ has the same distribution as 
$(S-B,S)$, where $S_t=\max\{B_s:0\leq s\leq t\}$ (to be specific, this holds with 
a particular choice of the normalization of $L$). We then see that Theorem \ref{subor-brown}
is a special case of Theorem  \ref{Levy-subor} and Proposition 
\ref{branchmecastable}. In that case, $\alpha=\beta=1/2$, and
we recover the formula $\psi(r)=c\,r^{3/2}$. 

\section*{Appendix: On the special Markov property}

In this appendix, we derive a more precise and more general form
of the special Markov property for the Brownian snake, which was first
stated in \cite{snakesolutions}. This result is closely related to
the special Markov property for superprocesses as stated by Dynkin
\cite[Theorem 1.6]{Dyn}, but the formulation in terms of the Brownian snake, although
less general,
gives additional information that is crucial for our purposes.

We consider the setting of \cite[Chapter V]{Zurich}. We let $\xi$ be a Markov process with
values in a Polish space $(E,d)$ with continuous sample paths. For every $x\in E$, the 
process $\xi$ starts from $x$ under the probability measure $P_x$. Analogously
to \eqref{contxi}, we assume that the following strong continuity assumption holds for
every $x\in E$ and $t\geq 0$, 
\begin{equation}
\label{contxi2}
E_x\Big[\Big(\sup_{r\leq t} d(x,\xi_r)\Big)^p\Big]\leq C\,t^{2+\ve},
\end{equation}
where $C>0$, $p>0$ and $\ve>0$ are constants.
According to \cite[Section IV.4]{Zurich}, this continuity assumption 
allows us to construct the Brownian snake $(W_s)_{s\geq 0}$ with continuous sample paths with
values in the space $\ww_E$ of all finite continuous paths in $E$
(the set $\ww_E$ is defined by the obvious generalization of the
beginning of Section \ref{sec:subord-maxi}, and we keep the notation
$\wh \w$ for the tip of a path $\w\in\ww_E$). The strong Markov property 
holds for $(W_s)_{s\geq 0}$, even without assuming that it holds for the
underlying spatial motion $\xi$. We again write 
$\P_x$ for the probability measure under which the Brownian snake starts
from
(the trivial path equal to) $x$, and
$\N_x$ for the excursion measure of the Brownian snake away from
$x$. It will also be useful to introduce conditional distributions of the
Brownian snake given its lifetime process. If $g:[0,\infty)\to[0,\infty)$
is a continuous function such that $g(0)=0$ and $g$ is locally H\"older with
exponent $\frac{1}{2}-\delta$ for every $\delta>0$, we write $\Q^{(g)}_x$ for the
conditional distribution under $\P_x$ of $(W_s)_{s\geq 0}$ knowing that
$\zeta_s=g(s)$ for every $s\geq 0$. See \cite[Chapter IV]{Zurich}, and note that
these conditional distributions are easily defined using the analog
in our general setting of property (b) stated at the beginning of 
Section \ref{sec:subord-maxi}. 

We now fix a connected open subset $D$ of $E$ and $x\in D$. We use the notation
$\tau=\inf\{t\geq 0:\xi_t\notin D\}$, and, for every $\w\in \ww_E$, $\tau(\w)=\inf\{t\geq 0:\w(t)\notin D\}$,
where in both cases $\inf\varnothing=\infty$. We assume that
$$P_x(\tau<\infty)>0,$$
and note that this implies that
$$\int_0^\infty \mathbf{1}_{\{\tau(W_s)<\zeta_s\}}\,\D s=\infty\,,\quad P_x\hbox{ a.s.}$$

We set, for every $s\geq 0$,
$$\eta_s:=\inf\Big\{t\geq 0: \int_0^t \mathbf{1}_{\{\zeta_r\leq \tau(W_r)\}}\,\D r >s\Big\}\;,\quad W^D_s:=W_{\eta_s}.$$
This definition makes sense $\P_x$ a.s. We let $\ff^D$ be the $\sigma$-field
generated by the process $(W^D_{s})_{s\geq 0}$ and the collection of all $\P_x$-negligible sets. Informally, 
$\ff^D$ represents the information provided by the paths 
$W_s$ before they exit $D$. 

\begin{lemma}
\label{MBreflected}
For every $s\geq 0$, set $\gamma_s=(\zeta_s-\tau(W_s))^+$, and 
$$\sigma_s=\inf\Big\{t\geq 0: \int_0^t \mathbf{1}_{\{\gamma_r>0\}}\,\D r \geq s\Big\}.$$
Under the probability measure $\P_x$, 
we have $\sigma_s<\infty$ for every $s\geq 0$, a.s., and the process $\Gamma_s:=\gamma_{\sigma_s}$
is distributed as a reflected Brownian motion in $\R_+$
and is independent of the $\sigma$-field $\ff^D$. 
\end{lemma}

This is essentially Lemma V.2 in \cite{Zurich}, except that the independence property
is not stated in that lemma. However a close look at the proof 
in \cite{Zurich} shows that the process $\Gamma_s$ is obtained as 
the limit of approximating processes (denoted by $\gamma_{\sigma^\ve_s}$ in \cite{Zurich})
which are independent of $\ff^D$ thanks to the strong Markov property of the Brownian snake.

We write $\ell^D(s)$ for the local time at $0$ of the process $\Gamma$, and define a process
with continuous nondecrasing sample paths by setting
$$L^D_s=\ell^D\Big(\int_0^s \mathbf{1}_{\{\gamma_r>0\}}\,\D r\Big).$$
Then (see \cite[Section V.1]{Zurich}),
\begin{equation}
\label{exit-local}
L^D_s=\lim_{\ve\to 0} \int_0^s \mathbf{1}_{\{\tau(W_r)<\zeta_r<\tau(W_r)+\ve\}}\,\D r,
\end{equation}
for every $s\geq 0$, $\P_x$ a.s. The process $(L^D_s)_{s\geq 0}$ is called the exit local time process from $D$. Notice that the measure $\D L^D_s$ is supported on
$\{s\geq 0: \tau(W_s)=\zeta_s\}$. 

We also set, for every $s\geq 0$,
$$\wt L^D_s=: L^D_{\eta_s}.$$

\begin{lemma}
\label{measurability-exit}
The random process $(\wt L^D_s)_{s\geq 0}$ is measurable with respect to the 
$\sigma$-field $\ff^D$.
\end{lemma}

This follows from the proof of \cite[Proposition 2.3]{snakesolutions} in the special case where $\xi$ is
$d$-dimensional Brownian motion. The argument however can be adapted to our
more general setting and we omit the details. 

Before stating the special Markov property, we need some additional notation.
For every $r\geq 0$, we set
$$\theta_r=\inf\{s\geq 0: \wt L^D_s>r \}.$$
We note that the process
$(\theta_r)_{r\geq 0}$ is $\ff^D$-measurable by Lemma \ref{measurability-exit}.

We now define the excursions of the Brownian snake outside $D$. 
We observe that, $\P_x$ a.e., the set
$$\{s\geq 0: \tau(W_s)<\zeta_s\}=\{s\geq 0: \gamma_s>0\}$$
is a countable union of disjoint open intervals, which we enumerate as
$(a_i,b_i)$, $i\in \N$. Here we can fix the enumeration by saying that we enumerate first
the excursion intervals with length at most $2^{-1}$ whose initial time is smaller than $2$, then
the excursion intervals with length at most $2^{-2}$ whose initial time is smaller than $2^2$ 
which have not yet been listed, and so on.  If we choose this enumeration procedure, the variables
$a_i,b_i$ are measurable with respect to the $\sigma$-field generated by $(\gamma_s)_{s\geq 0}$,
and the variables  $\int_0^{a_i} \mathbf{1}_{\{\gamma_r>0\}}\,\D r$ and $L^D_{a_i}$ are measurable with respect to the $\sigma$-field
generated by $(\Gamma_s)_{s\geq 0}$. 

From the properties of the Brownian 
snake, one has, $\P_x$ a.e. for every $i\in \N$ and every
$s\in[a_i,b_i]$, 
$$\tau(W_s)=\tau(W_{a_i})=\zeta_{a_i},$$
and more precisely all paths $W_s$, $s\in[a_i,b_i]$ coincide up to
their exit time from $D$. For every $i\in \N$, we then define an element
$W^{(i)}$ of the space of all continuous functions from $\R_+$ into $\ww_E$
by setting, for every $s\geq 0$,
$$W^{(i)}_s(t) := W_{(a_i+s)\wedge b_i}(\zeta_{a_i}+t),\quad \hbox{for }
0\leq t\leq \zeta^{(i)}_s:=\zeta_{(a_i+s)\wedge b_i}-\zeta_{a_i}.$$
By definition, the random variables $W^{(i)}$, $i\in \N$, are the excursions of the
Brownian snake outside $D$ (the word ``outside'' is a bit misleading here).
Notice that, for every $i\in\N$, $(\zeta^{(i)}_s)_{s\geq 0}$  is a measurable
function of $(\Gamma_s)_{s\geq 0}$. Indeed the processes  $(\zeta^{(i)}_s)_{s\geq 0}$, $i\in\N$
are just the excursions of $\Gamma$ away from $0$ enumerated as explained above.

\begin{theorem}
\label{SMPgeneral}
Under $\P_x$, conditionally on the $\sigma$-field $\ff^D$, the point measure
$$\sum_{i\in \N} \delta_{(L^D_{a_i}, W^{(i)})}(\D \ell,\D \omega)$$
is Poisson with intensity
$$\mathbf{1}_{[0,\infty)}(\ell)\,\D \ell \,\N_{\wh W^D_{\theta_\ell}}(\D \omega).$$
\end{theorem}

\proof
It is convenient to introduce the auxiliary Markov process defined by
$\xi^*_t=\xi_{t\wedge \tau}$. We observe that the Brownian snake 
associated with $\xi^*$
can be obtained by the formula
$$W^*_s(t)= W_s(t\wedge \tau(W_s))\,,\qquad 0\leq t\leq \zeta^*_s=\zeta_s.$$
Notice that $\gamma_s=(\zeta_s-\tau(W_s))^+ = (\zeta^*_s-\tau(W^*_s))^+$
is a measurable function of $W^*$, and recall that the intervals $(a_i,b_i)$ are
just the connected components of the complement of the zero set 
of $\gamma$. Consider then, independently for every $i\in \N$, 
a process $(\bar W^{(i)}_s)_{s\geq 0}$ which conditionally
on $W^*$ is distributed according to the probability measure
$$\Q_{\wh W^*_{a_i}}^{(\zeta^{(i)})},$$
where we recall our notation $\Q^{(g)}_x$ for the
conditional distribution under $\P_x$ of $(W_s)_{s\geq 0}$ knowing that
$\zeta_s=g(s)$ for every $s\geq 0$.
Define $\bar W_s$ for every $s\geq 0$ by setting $\bar W_s=W^*_s$
if $\gamma_s=0$ and, for every $i\in \N$, for every $s\in(a_i,b_i)$,
$$\bar W_s(t)=\left\{ 
\begin{array}{ll}
W^*_s(t)&\hbox{if }0\leq t\leq \tau(W^*_s),\\
\bar W^{(i)}_{s-a_i}(t-\tau(W^*_s))\quad&\hbox{if } \tau(W^*_s)\leq t\leq \zeta_s.
\end{array}
\right.
$$
A tedious but straightforward verification shows that the finite marginal distributions of the
process $(\bar W_s)_{s\geq 0}$ are the same as those of the process $(W_s)_{s\geq 0}$.
It follows that, conditionally on $(W^*_s)_{s\geq 0}$, the 
``excursions'' $W^{(i)}$ are independent and the conditional distribution
of $W^{(i)}$ is
$\Q_{\wh W^*_{a_i}}^{(\zeta^{(i)})}$. 

At this point, we claim that we have a.s. for every $i\in \N$,
\begin{equation}
\label{claim-tech}
W^*_{a_i}=W_{a_i}= W^D_{\theta_{L^D_{a_i}}}.
\end{equation}
The first equality in \eqref{claim-tech} is immediate. To get the second one, set
$$A_{a_i}:=\int_0^{a_i}\mathbf{1}_{\{\zeta_r\leq \tau(W_r)\}}\D r$$
to simplify notation.
We first note that
$$W_{a_i}= W_{b_i}=  W^D_{A_{a_i}},$$
because 
$\eta_{A_{a_i}}=b_i$
(the strong Markov property of the Brownian snake ensures that in each interval $[b_i,b_i+\ve]$, $\ve>0$, we can find a set of positive Lebesgue measure
of values of $s$ such that $\tau(W_s)=\infty$)
and we know that $W_{a_i}=W_{b_i}$. Thus our claim will follow if we can verify that
$$A_{a_i}=\theta_{L^D_{a_i}}.$$
On one hand the condition $s<A_{a_i}$ implies $\eta_s<a_i$
and $\wt L^D_s=L^D_{\eta_s}\leq L^D_{a_i}$. It follows that $A_{a_i}\leq \theta_{L^D_{a_i}}$.
On the other hand, suppose that $\theta_{L^D_{a_i}}> A_{a_i}$. We first note that
$$\wt L^D_{A_{a_i}}= L^D_{\eta_{A_{a_i}}}= L^D_{b_i}=L^D_{a_i},$$
where the last equality holds by the support property of $\D L^D_s$. 
Furthermore, the left limit of $r\mapsto \theta_r$ at $\wt L^D_{A_{a_i}}$ is smaller than or
equal to $A_{a_i}$ by construction. So the condition $\theta_{L^D_{a_i}}> A_{a_i}$
means that $L^D_{a_i}=\wt L^D_{A_{a_i}}$ is a discontinuity point of
$r\mapsto \theta_r$. However, we noticed that the random variables $L^D_{a_i}$
are measurable functions of $(\Gamma_s)_{s\geq 0}$ 
and therefore independent of $\ff^D$. Since $(\theta_s)_{s\geq 0}$ is 
measurable with respect to $\ff^D$, and since $s\mapsto \theta_r$
has only countably many discontinuity times, the (easy) fact that the law
of $L^D_{a_i}$ has no atoms implies that, with $\P_x$-probability one, 
$L^D_{a_i}$ cannot be a discontinuity time of
$r\mapsto \theta_r$. This contradiction completes the proof of our claim.

Next, let $U$ be a bounded $\ff^D$-measurable real random variable, and let
$g$ and $G$ be nonnegative random variables defined respectively 
on $\R_+$ and on the space of continuous functions from $\R_+$
into $\ww_E$. By conditioning first with respect to $W^*$, we get
\begin{equation}
\label{SMP1}
\E_x\Big[ U\times \exp\Big(-\sum_{i\in\N} g(L^D_{a_i}) G(W^{(i)})\Big)\Big]
= \E_x\Bigg[ U\times \prod_{i\in \N} \Q_{\wh W_{a_i}}^{(\zeta^{(i)})}\Big(e^{-g(L^D_{a_i}) G(\cdot)}
\Big)\Bigg],
\end{equation}
noting that an $\ff^D$-measurable real variable coincides $\P_x$ a.s. with a function 
of $W^*$. Using \eqref{claim-tech}, we see that, for every $i\in \N$, 
$$\Q_{\wh W_{a_i}}^{(\zeta^{(i)})}\Big(e^{-g(L^D_{a_i}) G(\cdot)}
\Big)=
\Q_{\wh W^D_{\theta_{L^D_{a_i}}}}^{(\zeta^{(i)})}\Big(e^{-g(L^D_{a_i}) G(\cdot)}
\Big)$$
is a measurable function (that does not depend on $i$) of the pair $(L^D_{a_i},
(\zeta^{(i)}_s)_{s\geq 0})$
and of the process $(W^D_{\theta_r})_{r\geq 0}$, which is $\ff^D$-measurable.
The point measure 
$$\sum_{i\in \N} \delta_{(L^D_{a_i},(\zeta^{(i)}_s)_{s\geq 0})},$$
which is just the point measure of excursions of the reflected Brownian motion $\Gamma$,
is Poisson with intensity $\D \ell\,\bn(\D e)$, where $\bn(\D e)$ is as previously the 
It\^o excursion measure. Since the latter point measure is
independent of $\ff^D$ (by Lemma \ref{MBreflected}), we can now 
condition with respect to $\ff^D$, applying 
the exponential formula for Poisson measures, to get that the
quantities in \eqref{SMP1} are equal to
$$\E_x\Bigg[ U\times \exp\Big(-\int_0^\infty \D \ell \, \N_{\wh W^D_{\theta_\ell}}\Big(1-
e^{-g(\ell)G(\cdot)}\Big)\Big)\Bigg].$$
The statement of the theorem follows. \endproof

In the preceding sections, we use a version of Theorem
\ref{SMPgeneral} under the excursion measure $\N_x$,
which we will now state as a corollary. We set
$$T_D=\inf\{t\geq 0: \tau(W_s)<\infty\},$$
and observe that
$$0<\N_x(T_D<\infty)<\infty$$
(if the quantity $\N_x(T_D<\infty)$ were infinite, excursion theory would give a contradiction with
the fact that the Brownian snake has continuous paths under $\P_x$). Then we note
that the conditional probability measure $\N_x^D:=\N_x(\cdot\mid T_D<\infty)$
can be interpreted as the law under $\P_x$ of the first Brownian excursion
away from $x$ that exits $D$. Thanks to this observation, we can make sense
of the exit local time process $(L^D_s)_{s\geq 0}$ under $\N^D_x$ by formula 
\eqref{exit-local}. The definition of $(\eta_s)_{s\geq 0}$
and $(W^D_s)_{s\geq 0}$ remains the same, 
and we can set $\wt L^D_s=L^D_{\eta_s}$ as previously. The difference of course
is the fact that $\wt L^D_\infty=L^D_\infty= L^D_\sigma$ is now finite
$\N^D_x$ a.s. So we may define $\theta_s=\inf\{r>0: \wt L^D_r>s\}$ 
only for $r<L^D_\sigma$. 

After these observations, we may state the version of the special Markov property
under $\N^D_x$. We slightly abuse notation by still writing 
$(W^{(i)})_{i\in\N}$ and $(a_i,b_i)_{i\in \N}$ for the excursions outside $D$
and the associated intervals, which are defined in exactly
the same way as previously. The $\sigma$-field $\ff^D$ is again
the $\sigma$-field generated by $(W^D_s)_{s\geq 0}$ but should be 
completed here with the $\N_x$-negligible sets.

\begin{corollary}
\label{SMPexcursion}
Under $\N^D_x$, conditionally on the $\sigma$-field $\ff^D$, the point measure
$$\sum_{i\in \N} \delta_{(L^D_{a_i}, W^{(i)})}(\D \ell,\D \omega)$$
is Poisson with intensity
$$\mathbf{1}_{[0,L^D_\sigma]}(\ell)\,\D \ell \,\N_{\wh W^D_{\theta_\ell}}(\D \omega).$$
\end{corollary}

Note that $L^D_\sigma=\wt L^D_\infty$ is $\ff^D$-measurable. 
The corollary is a straightforward consequence of Theorem \ref{SMPgeneral}, 
interpreting $\N^D_x$ as the law under $\P_x$ of the first Brownian excursion
away from $x$ that exits $D$. We
omit the details. 

If we are only interested in the point measure $\sum_{i\in I}\delta_{W^{(i)}}$, we can state
the preceding corollary in a slightly different form, by introducing the notion of
the exit measure: The exit measure from $D$ is the random measure $\zz^D$
on $\partial D$ defined (under $\N_x$ or under $\N^D_x$) by the formula
$$\langle \zz^D,g\rangle =\int_0^\infty \D L^D_s\, g(\wh W_s).$$
See \cite[Section V.1]{Zurich}. Note that the total mass of $\zz^D$
is $L^D_\sigma$. A change of variables shows that we have as well
$$\langle \zz^D,g\rangle = \int_0^\infty \D \wt L^D_s\,g(\wh W^D_s)= \int_0^{L^D_\sigma} \D \ell \,g(\wh W^D_{\theta_\ell}).$$
In particular, $\zz^D$ is $\ff^D$-measurable by Lemma \ref{measurability-exit}. The following result
is now an immediate consequence of Corollary \ref{SMPexcursion}.

\begin{corollary}
\label{SMPusual}
Under $\N^D_x$, conditionally on the $\sigma$-field $\ff^D$, the point measure
$$\sum_{i\in \N} \delta_{ W^{(i)}}(\D \omega)$$
is Poisson with intensity
$$\int \zz^D(\D y)\,\N_y(\D \omega).$$
\end{corollary}

This is the form of the special Markov property that appears in \cite{snakesolutions}.

\end{document}